\documentclass[11pt]{article}
\usepackage[latin1]{inputenc}
\usepackage{epsfig}
\usepackage{color}
\usepackage[british,english]{babel}
\usepackage{amsthm}
\usepackage{amsmath}
\usepackage{amsfonts}
\usepackage{amssymb}
\usepackage{graphicx}
\newtheorem{corollary}{Corollary}[section]
\newtheorem{definition}[corollary]{Definition}

\newtheorem{lemma}[corollary]{Lemma}
\newtheorem{proposition}[corollary]{Proposition}
\newtheorem{remark}[corollary]{Remark}
\newtheorem{theorem}[corollary]{Theorem}
\newfont{\sBlackboard}{msbm10 scaled 900}

\newcommand{\mylabel}[1]{\label{#1}
            \ifx\undefined\stillediting
            \else \fbox{$#1$}\fi }
\newcommand{\BE}{\begin{equation}}

\newcommand{\EEQ}{\end{equation}}
\newcommand{\rfb}[1]{\mbox{\rm
   (\ref{#1})}\ifx\undefined\stillediting\else:\fbox{$#1$}\fi}

\newfont{\Blackboard}{msbm10 scaled 1200}

\newfont{\roma}{cmr10 scaled 1200}

\def\CC{\rm \hbox{C\kern-.56em\raise.4ex
         \hbox{$\scriptscriptstyle |$}\kern+0.5 em }}

\newcommand{\ep}{\varepsilon}

\def\n{|\kern -.05cm{|}\kern -.05cm{|}}

\newcommand{\mm}    {{\hbox{\hskip 0.5pt}}}

\newcommand{\bluff} {{\hbox{\raise 15pt \hbox{\mm}}}}
\usepackage{fancyhdr}

\makeatletter
\def\section{\@startsection {section}{1}{\z@}{-3.5ex plus -1ex minus
    -.2ex}{2.3ex plus .2ex}{\large\bf}}
\makeatother
\def\be{\begin{equation}}
\def\ee{\end{equation}}

\date{ }
\begin{document}
%%-----------------------------
%%      the top matter
%%-----------------------------
\thispagestyle{empty}
\title{\bf Modeling of a non-Newtonian thin  film\\
 passing a thin porous medium}
\maketitle
%non-zero boundary conditions for microrotations}\maketitle
\vspace{-50pt}
\author{ \center  Mar\'ia ANGUIANO\\
Departamento de An\'alisis Matem\'atico. Facultad de Matem\'aticas.\\
Universidad de Sevilla, 41012-Sevilla (Spain)\\
anguiano@us.es\\}
%\author{ \center  Matthieu BONNIVARD\\
%Univ Lyon, Ecole Centrale de Lyon, CNRS, Univ Claude Bernard Lyon 1, \\
%Univ Jean Monnet, INSA Lyon, Institut Camille Jordan, UMR5208, 69130 Ecully (France)\\
%matthieu.bonnivard@ec-lyon.fr\\}
\author{ \center  Francisco Javier SU\'AREZ-GRAU\\
Departamento de Ecuaciones Diferenciales y An\'alisis Num\'erico. Facultad de Matem\'aticas.\\
Universidad de Sevilla,  41012-Sevilla (Spain)\\
fjsgrau@us.es\\}

\vskip20pt
\begin{abstract} This theoretical study deals with asymptotic behavior of a coupling between a thin film of fluid and an adjacent thin porous medium. We assume that the size of the microstructure of the porous medium is given by a small parameter $0<\ep\ll 1$,  the thickness of the thin porous medium is defined by a parameter $0<h_\ep\ll 1$, and the thickness of the thin film is defined by a small parameter $0<\eta_\ep\ll 1$, where $h_\ep$ and $\eta_\ep$ are devoted to tend to zero when $\ep\to 0$. In this paper,  we consider the case of a non-Newtonian fluid governed by the incompressible Stokes equations with power law viscosity of flow index $r\in (1, +\infty)$, and we prove that there exists a critical regime, which depends on $r$, between  $\ep$,   $\eta_\ep$ and   $h_\ep$. More precisely,  in this critical regime given by $h_\ep\approx \eta_\ep^{2r-1\over r-1}\ep^{-{r\over r-1}}$, we prove that the effective flow when $\ep\to 0$ is described by a 1D Darcy law coupled with a 1D Reynolds law.  \end{abstract}
%
%\begin{resume} ... \end{resume}
%
\bigskip\noindent
 {\small \bf AMS classification numbers:}  74Q10, 76A05, 76A20, 76S05. 

%%-----------------------------
%%      your text
%%-----------------------------
\section {Introduction}\label{S1}

% IMPORTANCE OF THE STUDY IN APPLICATIONS

In this paper, we consider a incompressible viscous 2D non-Newtonian fluid in a domain $D_\ep$ composed by two parts in contact: a periodic thin porous medium $\Omega_\ep$ with characteristic size of the pores $0<\ep\ll1$ and thickness of the domain $0<h_\ep\ll 1$, and a thin film $I_\ep$ with thickness $0<\eta_\ep\ll 1$, where $h_\ep$ and $\eta_\ep$ are devoted to tend to zero when $\ep\to 0$ (see Figure \ref{fig:domainD} for more details). Drilling and hydraulic fracturing fluids used in the oil industry are usually non-Newtonian liquids. Therefore during well drilling or hydraulic fracturing operations, the non-Newtonian drilling muds or hydraulic fluids will
infiltrate into permeable formations surrounding the wellbore, which may seriously damage the formation. The rheological behavior of drilling muds, cement slurries and hydraulic fracturing fluids is often described by a power-law model (see Cloud and Clark \cite{Cloud}, Shah \cite{Shah}). The importance of modeling flow
of non-Newtonian fluids from the wellbore into the surrounding formations has been recognized in the
industry.

One way to study this problem is to use the homogenization theory, which has been applied to
the study of perforated materials for a long time (see for instance classical studies of Allaire \cite{Allaire0}, Sanchez-Palencia \cite{Sanchez} and  Tartar \cite{Tartar} in the case of Newtonian fluids, or  Bourgeat {\it et al.} \cite{BourgeatGipoloux} and Bourgeat and Mikeli\'c \cite{Bourgeat1} in the case of non-Newtonian fluids with viscosity given by the power law or the Carreau law). The question of a porous medium in contact with a thin film
with properties different from those of the rest of the material has been the subject of many studies
previously.

% ANTECEDENTES

Let us make a recollection of some previous results in relation to the objective of this paper.  Bourgeat {\it et al.} \cite{BourgeatPalokaMikelic} considered the asymptotic behavior of the solution of the 2D Newtonian Stokes system in a porous medium with thickness of order one, with characteristic size of the pores $\ep$ and containing a thin fissure of thickness $\eta_\ep$, with $\eta_\ep$ devoted to tend to zero with $\ep$. It was proved the existence of a critical regime  given by 
\begin{equation}\label{critical_fissure_Newtonian}
\eta_\ep\approx \ep^{2\over 3},
\end{equation} where the coupling effect appears and the effective flow behaves like a 2D Darcy flow in the porous medium coupled with a 1D Reynolds problem.  We  refer to Bourgeat {\it et al.} \cite{BourgeatElAmriTapiero} for preliminary results on this problem, and to Zhao and Yao \cite{ZhaoYao1, ZhaoYao2} for the extension of this result to the non-stationary Stokes case and the Navier-Stokes  case, respectively. We also refer to Bourgeat \cite{BourgeatPalokaMikelicLaplace} and Bourgeat and Tapi\'ero \cite{BourgeatTapieroLaplace} for a similar problem but for the Laplace equation.   

Moreover, this problem was also generalized in Anguiano and Su\'arez-Grau \cite{Anguiano_CMS} to the case of a non-Newtonian Stokes flow with viscosity given by the power law with flox index $r$ satisfying $1<r<+\infty$, where the critical regime  is now given by 
\begin{equation}\label{critical_fissure_nonNewtonian}
\eta_\ep\approx \ep^{r\over 2r-1},
\end{equation}
which agrees with (\ref{critical_fissure_Newtonian}) for $r=2$.  In this case, the effective flow behaves like a 2D nonlinear Darcy flow in the porous medium coupled with a 1D nonlinear Reynolds problem.  We also refer to Anguiano \cite{Anguiano_fis2} for the case of a non-stationary non-Newtonian flow in a porous medium containing a thin fissure, where in the critical regime (\ref{critical_fissure_nonNewtonian}), the flow is described by a time-dependent non-linear Reynolds problem coupling the effects of the porous medium with those of the free part.

%%% THIN POROUS MEDIUM

On the other hand, the derivation of effective laws for fluids in porous domains with small thickness (the so-called {\it thin porous medium}) is attracting much attention, see for instance Almqvist {\it et al.} \cite{Almqvist}, Anguiano \cite{Anguiano2}, Anguiano and Bunoiu \cite{Ang-Bun2}, Anguiano {\it et al.} \cite{Anguiano_Bonnivard_SG, Anguiano_Bonnivard_SG2}, Anguiano and Su\'arez-Grau \cite{Anguiano_SuarezGrau, Anguiano_SuarezGrau2, Anguiano_SG_NHM, Anguiano_SG_Lower, Anguiano_SG_sharp}, Fabricius {\it et al.} \cite{Fabricius, Fabricius0}, Fabricius and Gahn \cite{Fabricius3}, Forslund {\it et al.} \cite{Forslund},  Jouybari and  Lundstr$\ddot{\rm o}$m \cite{Jouybari}, Mei and Vernescu \cite{Mei}, Su\'arez-Grau \cite{SG_MANA, SuarezGrau1} or Zhengan and  Hongxing \cite{Zhengan}.  A thin porous medium can be defined as a bounded domain confined between two parallel plates with a distance $h_\ep$, perforated by periodically distributed obstacles of size $\ep$, with $h_\ep$ devoted to tend to zero when $\ep\to 0$.

In this context, the modeling of a Newtonian flow in a thin porous medium and an adjacent thin film flow described by Figure \ref{fig:domainD}, which is the domain we are interested in this paper,  was considered in Bayada {\it et al} \cite{BayadaThinThin}. Thus,  considering  three positive and small parameters $\ep$, $h_\ep$ and $\eta_\ep$ (where $h_\ep$ and $\eta_\ep$ are devoted to tend to zero), where $\ep$ is the size of the microstructure,  $h_\ep$ is the thickness of the thin porous medium and $\eta_\ep$ is the thickness of the thin film, it was proved  the existence of a critical regime between these parameters given by
\begin{equation}\label{Critical_thin_thin_Newtonian}
h_\ep\approx  \eta_\ep^3 \ep^{-2},
\end{equation}
and  an effective modified Reynolds equation (a 1D Darcy problem coupled with a 1D Reynolds problem) was derived. Observe that if the thickness of the porous medium $h_\ep\equiv 1$, then the critical regime (\ref{Critical_thin_thin_Newtonian}) coincides with that critical one given in  (\ref{critical_fissure_Newtonian}). We also refer to Anguiano and Su\'arez-Grau \cite{Anguiano2} for the derivation of a coupled Darcy--Reynolds equation for a fluid flow in a thin porous medium including a fissure, where the microstructure of thin porous medium is a collection of small cylinders (see Anguiano \cite{Anguiano_fis} for the non-stationary case).

% OBJETIVE OF THE PAPER

The goal of this paper is to generalize the result described in  \cite{BayadaThinThin} to the case of a non-Newtonian fluid with a viscosity described by the power law with flow index $r\in (1,+\infty)$, which is important for industrial applications as described above. We prove that there exists a critical regime between these parameters given by 
\begin{equation}\label{Critical_thin_thin_nonNewtonian}
h_\ep\approx  \eta_\ep^{2r-1\over r-1} \ep^{-{r\over r-1}},
\end{equation}
and we derive an effective nonlinear limit problem, i.e. a modified nonlinear Reynolds problem coupling the effects of the thin porous medium (1D nonlinear Darcy problem) and the thin film (1D nonlinear Reynolds problem) for the limit pressure (see Theorem \ref{main_thm} for more details). We observe that if $r=2$, then the critical regime (\ref{Critical_thin_thin_nonNewtonian}) coincides with the critical one given in  (\ref{Critical_thin_thin_Newtonian}). Also,  if the thickness of the porous medium $h_\ep\equiv 1$, then the critical regime (\ref{Critical_thin_thin_nonNewtonian}) coincides with the critical one given in  (\ref{critical_fissure_nonNewtonian}).  To prove this result, we first derive global {\it a priori} estimates of the velocity and pressure and also, particular {\it a priori} estimates in both media, which let us find the critical regime (\ref{Critical_thin_thin_nonNewtonian}). Then, in this critical regime, we study independently the asymptotic behavior in the thin porous medium and in the thin film. Finally, we  deduce that the pressure is continuous in the interface and  derive the coupled effective problem for the limit pressure. We have introduced the following novelties in this work  with respect to previous studies to study the asymptotic behavior in the thin porous medium:  a new restriction operator $\mathcal{R}^\ep_r$ to extend the pressure in the thin porous medium to the thin domain without microstructure (see Subsection  \ref{sec:restriction}), and a new version of the unfolding method  (for classical versions see  Cioranescu {\it et al.} \cite{Ciora2,Cioran-book}) to capture the effects of the microstructure of the thin porous medium (see Subsection \ref{sec:unfolding}). All this is combined with dimension reduction techniques and monotonicity arguments to be able to pass to the limit when $\ep\to 0$ and so, to derive the modified Reynolds equation.

We think that this theoretical study provides a quite complete description of the asymptotic behavior of generalized Newtonian fluids with power law viscosity through a thin film passing a thin porous medium, which provides a model  amenable for the numerical simulations. For this reason,  we believe that it could also prove useful in the engineering practice as well.

Finally, we comment the structure of the paper. In Section \ref{sec:domain} we introduce the domain considered and the statement of the problem. In Section \ref{sec:estimates},  we derive  {\it a priori} estimates of the velocity and pressure. Thanks to the estimates, we find the critical regime (\ref{Critical_thin_thin_nonNewtonian}).  Assuming the critical regime, in Section \ref{sec:thinporousmedium} we study the asymptotic behavior of the problem in the thin porous medium (extending the pressure by duality arguments using the restriction operator $\mathcal{R}^\ep_r$ and using the version of the unfolding method to capture the effects of the microstructure), and in Section \ref{sec:thinfilm}, we study the asymptotic behavior of the problem in the thin film (which differs a little bit from the classical study of the asymptotic behavior of a non-Newtonian fluid in a thin domain). Finally, in Section \ref{sec:mainresult} we deduce that the pressure is continuous in the interface of both media and we derive the modified Reynolds problem coupling the effects of the thin porous medium and the thin film, which is given Theorem \ref{main_thm}. We finish the paper with a conclusion section and a list of references.

\section{Formulation of the problem and preliminaries} \label{sec:domain} 
\subsection{Geometrical setting}

\begin{figure}[h!]
\begin{center}
\includegraphics[width=8cm]{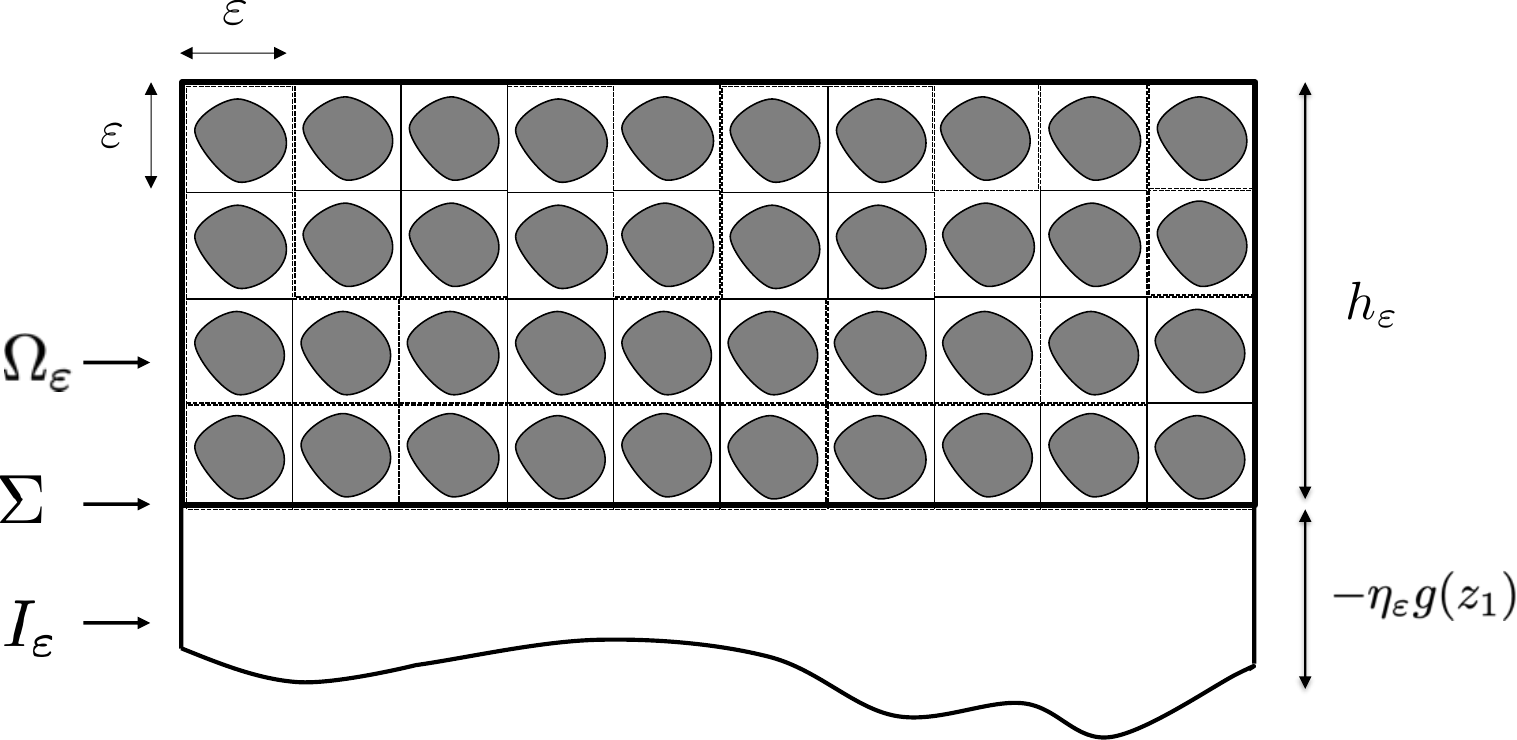}
\end{center}
\vspace{-0.4cm}
\caption{View of the domain $D_\ep$.}
\label{fig:domainD}
\end{figure}

\begin{figure}[h!]
\begin{center}
\includegraphics[width=8cm]{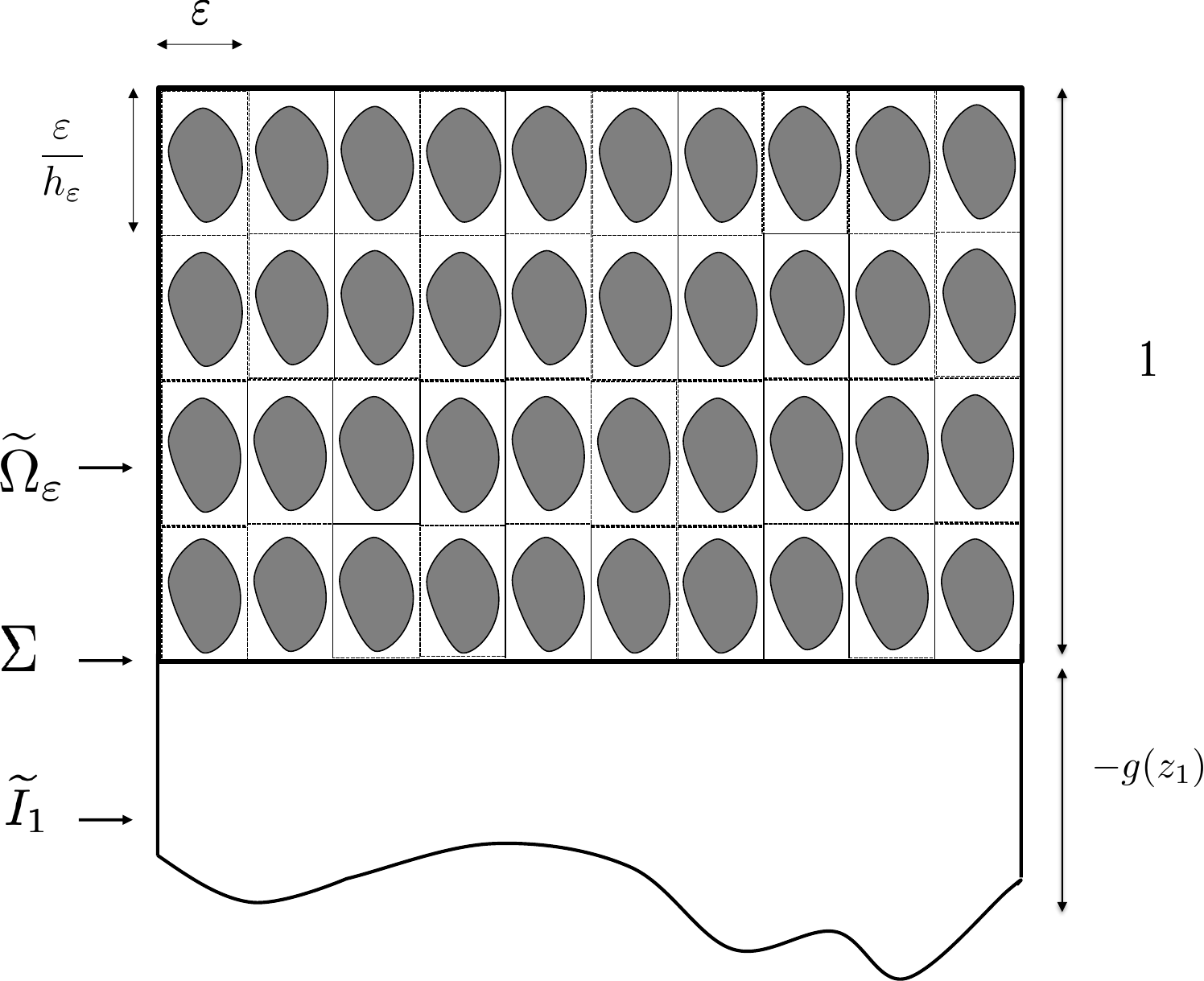}
\end{center}
\caption{View of the rescaled domain $\widetilde D_{\varepsilon}$.}
\label{fig:domain2}
\end{figure}

Let $\omega\in (-1/2,1/2)\subset \mathbb{R}$. We consider three positive and small parameters $\ep$, $h_\ep$ and $\eta_\ep$ (where $h_\ep$ and $\eta_\ep$ are devoted to tend to zero when $\ep\to0$) satisfying the following relation
\begin{equation}\label{parameters}
\lim_{\ep\to 0} {\ep\over h_\ep}=0,\quad \left(\hbox{i.e.}\quad \ep\ll h_\ep\right),\quad \hbox{and}\quad \lim_{\ep\to 0} {\ep\over \eta_\ep}=0,\quad \left(\hbox{i.e.}\quad \ep\ll \eta_\ep\right).
\end{equation}
We consider $D_\ep\subset \mathbb{R}^2$ to be an open set of the following form
$$D_\ep=\Omega_\ep \cup \Sigma \cup I_\ep,$$
where $\Omega_\ep$ is a thin porous medium and $I_\ep$ is a thin layer without obstacles (see Figure \ref{fig:domainD}). Moreover,  $\Sigma$ is the interface between the thin porous medium and the thin film and is defined by 
$$\Sigma=\omega\times \{x_2=0\}.$$
Below, we describe  subdomains $\Omega_\ep$ and $I_\ep$:
\begin{itemize}
\item To describe the thin porous medium $\Omega_\ep$, we consider the parameters $\ep$ and $h_\ep$ satisfying (\ref{parameters}). We consider a thin layer of height $h_\ep$ which is perforated by $\ep$-periodic distributed obstacles of size $\ep$. The thin layer without microstructure is denoted by $Q_\ep$, i.e.
\begin{equation}\label{Qep}
Q_\ep=\omega\times (0,h_\ep).
\end{equation} 
Let us now give a better description of the microstructure of the thin layer. We denote $Y=(-1/2,1/2)^2$ the unitary cube in $\mathbb{R}^2$ as the reference cell and  $T$ an open connected subset of $Y$ with a smooth boundary  $\partial T$ such that $\overline T\subset Y$. We denote $Y_f=Y\setminus \overline T$. Thus, for $k\in\mathbb{Z}^2$, each cell $Y_{k,\varepsilon}=\varepsilon k+\varepsilon Y$ is similar to the unit cell $Y$ rescaled to size $\varepsilon$ and $T_{k,\varepsilon}=\varepsilon k+\varepsilon T$ is similar to $T$ rescaled to size $\varepsilon$. We denote $Y_{f_k,\varepsilon}=Y_{k,\varepsilon}\setminus \overline T_{k,\varepsilon}$ (see Figure \ref{fig:cell}).

\begin{figure}[h!]
\begin{center}
\includegraphics[width=7cm]{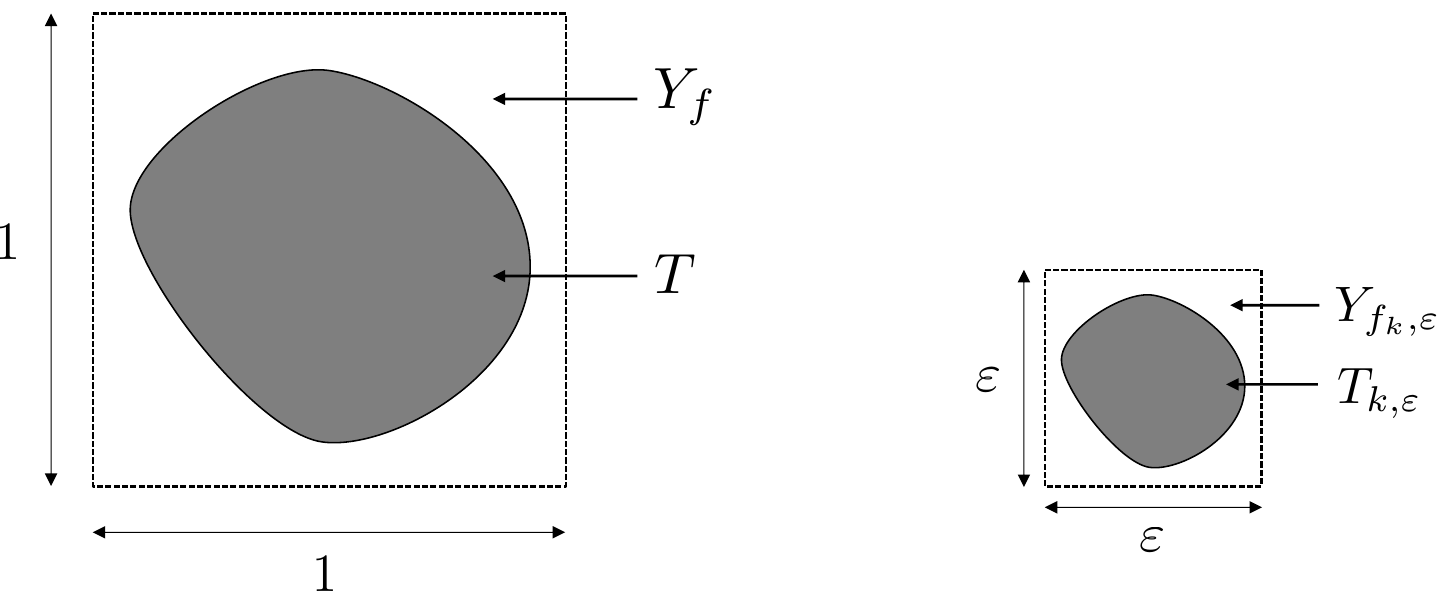}
\end{center}
\vspace{-0.4cm}
\caption{View of the reference cell  $Y$ (left) and the rescaled cell $Y_{k,\ep}$ (right).}
\label{fig:cell}
\end{figure}
We denote by $\tau(\overline T_{k,\ep})$ the set of all translated images of $\overline T_{k,\ep}$, i.e. the set $\tau(\overline T_{k,\ep})$ represents the obstacles in $\mathbb{R}^2$. 

\noindent The thin porous media $\Omega_\ep$ is defined by  (see Figure \ref{fig:domainD})
\begin{equation}\label{Omegaep}
\Omega_\ep=Q_\ep\setminus \bigcup_{k\in \mathcal{K}_\ep}\overline T_{k,\varepsilon},
\end{equation} where   $\mathcal{K}_\ep:=\left\{k\in\mathbb{Z}^3\,:\, Y_{k,\varepsilon}\cap Q_\ep\neq \emptyset\right\}$. By construction,  $\Omega_\varepsilon$ is a periodically perforated channel with obstacles of the same size as the period.
We make the assumption that the obstacles $\tau(\overline T_{k,\ep})$ do no intersect the boundary $\partial Q_\ep$. We denote by $T_{\varepsilon}$ the set of all the obstacles contained in $\Omega_\ep$. Then,   $T_\ep$ is a finite union of obstacles, i.e. 
 $$T_{\varepsilon}=\bigcup_{k\in \mathcal{K}_\ep}\overline T_{k,\ep}.$$
\noindent As usual when we deal with thin domains, we will use the dilatation in the variable $x_2$ given by
\begin{equation}\label{dilatacion_p}
z_1= x_1,\quad z_2={x_2\over h_\ep},\quad \forall\,x\in \Omega_\ep.
\end{equation}
Then, we define the rescaled porous media  $\widetilde \Omega_\ep$ by (see Figure \ref{fig:domain2})
\begin{equation}\label{Omega_tilde}
\widetilde\Omega_\ep=\left\{z=(z_1,z_2)\in \mathbb{R}^2\,:\, (z_1,h_\ep z_2)\in \Omega_\ep\right\}.
\end{equation}
We also introduce the rescaled sets $\widetilde Y_{k,\ep}$ by  (see Figure \ref{fig:domain2})
$$\widetilde Y_{{k},\ep}=\left\{z\in\mathbb{R}^2\,:\, (z_1,\ep z_2)\in Y_{k,\varepsilon}\right\},$$
and, in the same way, we define the rescaled fluid part $\widetilde Y_{f_{k},\ep}$, the rescaled solid  part $\widetilde T_{k,\varepsilon}$ of $\widetilde Y_{{k},\ep}$ and the union of rescaled obstacles $\widetilde T_\ep$.
Finally,  we denote by $\Omega$ the domain with fixed height without microstructure, i.e.
$$
\Omega=\omega\times (0,1). 
$$

\item To describe the thin layer $I_\ep$, we consider the positive and small parameter $\eta_\ep$ satisfying (\ref{parameters}). We define $I_\ep$ as follows
\begin{equation}\label{SetIep}
I_\ep=\left\{x=(x_1, x_2)\in \mathbb{R}^2\,:\, x_1\in \omega,\ g_\ep(x_1)<x_2<0\right\},
\end{equation}
where the function $g_\ep$ is given by
$$ g_\ep(x_1)=-\eta_\ep g(x_1),\quad \forall\, x_1\in \omega.$$
We define the lower boundary by 
$$\Gamma_g^\ep=\left\{x=(x_1, x_2)\in \mathbb{R}^2\,:\, x_1\in \omega,\ x_2=g_\ep(x_1)\right\}.
$$
Moreover, the following assumptions concerning the function $g$ are made:
\begin{equation}\label{Funcg}
g\in C(\omega),\quad 0< a\leq g(x_1)\leq b,\quad \forall x_1\in \omega\quad (\hbox{with }a,b>0).
\end{equation}
\noindent As before, to rescale $I_\ep$ in a set with fixed thickness, we will use the dilatation in the variable $x_2$ given by
\begin{equation}\label{dilatacion_t}
z_1= x_1,\quad z_2={x_2\over \eta_\ep},\quad  \forall\,x\in I_\ep.
\end{equation}
Then, we define the rescaled  domain $\widetilde I_1$ by 
\begin{equation}\label{SetI1}
\widetilde I_1=\left\{z=(z_1,z_2)\in \mathbb{R}^2\,:\, z_1\in\omega,\ -g(z_1)<z_2<0\right\},
\end{equation}
and  the  rescaled lower boundary by 
$$\Gamma_g=\left\{z=(z_1, z_2)\in \mathbb{R}^2\,:\, z_1\in \omega,\ z_2=-g(z_1)\right\}.
$$
\end{itemize}
Finally, we define the domain with microstructure by
$$
\Lambda_\ep=Q_\ep\cup \Sigma\cup I_\ep,
$$
the rescaled domain with microstructure (see Figure \ref{fig:domain2}) by
$$\widetilde D_\ep=\widetilde\Omega_\ep\cup \Sigma \cup \widetilde I_1,$$
and the whole rescaled domain without microsturcture   by 
$$D=\Omega\cup \Sigma \cup \widetilde I_1.
$$

\subsection{Some notation} Let us consider a vectorial function $\varphi=(\varphi_1, \varphi_2)$ and a scalar function $\psi$.  We have denoted by  
$\mathbb{D}:\mathbb{R}^2\to \mathbb{R}^2_{\rm sym}$   the symmetric part of the velocity gradient, that is
$$\mathbb{D}[\varphi]={1\over 2}(D\varphi+(D\varphi)^T)=\left(\begin{array}{cc}
\partial_{x_1}\varphi_1 &   {1\over 2}(\partial_{x_1}\varphi_2 + \partial_{x_2}\varphi_1)\\
\\
 {1\over 2}(\partial_{x_1}\varphi_2 + \partial_{x_2}\varphi_1) & \partial_{x_2}\varphi_2
\end{array}\right),$$
and also,  we have used the following operators 
$$\Delta\varphi=\partial_{x_1}^2\varphi+\partial_{x_2}^2\varphi,\quad {\rm div}(\varphi)=\partial_{x_1}\varphi_1+\partial_{x_2}\varphi_2,\quad \nabla\psi=(\partial_{x_1}\psi, \partial_{x_2}\psi)^t. $$
 For a vectorial function $\widetilde \varphi=(\widetilde \varphi_1,\widetilde \varphi_2)$ and a scalar function $\widetilde \psi$ obtained respectively from $\varphi$ and $\psi$ by using the rescaling (\ref{dilatacion_p}) in the set $\Omega_\ep$, we will denote  
$$\partial_{z_1}[\widetilde\varphi]={1\over 2}\left(\partial_{z_1}\widetilde\varphi+(\partial_{z_1}\widetilde\varphi)^T\right)
=\left(\begin{array}{cc} \partial_{z_1}\widetilde\varphi_1 & {1\over 2}\partial_{z_1}\widetilde\varphi_2\\
\\
{1\over 2}\partial_{z_1}\widetilde\varphi_2 & 0
\end{array}\right),\quad \partial_{z_2}[\widetilde\varphi]=\left(\begin{array}{cc} 0 & {1\over 2}\partial_{z_2}\widetilde \varphi_1\\
\\
{1\over 2}\partial_{z_2}\widetilde\varphi_1 & \partial_{z_2}\widetilde\varphi_2\end{array}\right),$$
and then, 
\begin{equation}\label{defD}\mathbb{D}_{h_\ep}[\widetilde\varphi]=\partial_{z_1}[\widetilde\varphi] + h_\ep^{-1}\partial_{z_2}[\widetilde\varphi]=\left(\begin{array}{cc} \partial_{z_1} \widetilde\varphi_1 &   {1\over 2}(\partial_{z_1}\widetilde\varphi_2 + h_\ep^{-1}\partial_{z_2}\widetilde\varphi_1)\\
\\
 {1\over 2}(\partial_{z_1}\widetilde\varphi_2 + h_\ep^{-1}\partial_{z_2}\widetilde\varphi_1) &  h_\ep^{-1}\partial_{z_2}\widetilde\varphi_2
\end{array}\right).
\end{equation}
Also, we define the  operators $\Delta_{h_\ep}$, $D_{h_\ep}$, $\nabla_{h_\ep}$ and ${\rm div}_{h_\ep}$ as follows
$$\begin{array}{c}
\displaystyle (D_{ h_\ep}\widetilde \varphi)_{i,1}=\partial_{x_1}\widetilde \varphi_i,\quad \displaystyle (D_{ h_\ep}\widetilde \varphi)_{i,2}=h_\ep^{-1}\partial_{z_2}\widetilde \varphi_i\hbox{ for }i=1,2,\\
\\
\displaystyle \Delta_{h_\ep}\widetilde\varphi=\partial_{z_1}^2\widetilde\varphi+h_\ep^{-2}\partial_{z_2}^2\widetilde\varphi,\quad \displaystyle  {\rm div}_{h_\ep} (\widetilde\varphi)=\partial_{z_1}\widetilde\varphi_1+h_\ep^{-1}\partial_{z_2}\widetilde\varphi_2,\\
\\\displaystyle
\nabla_{ h_\ep}\widetilde\psi=(\partial_{x_1}\widetilde\psi,h_\ep^{-1}\partial_{z_2}\widetilde\psi)^t.
\end{array}$$
Similarly, we define the operators $\mathbb{D}_{\eta_\ep}$, $\Delta_{\eta_\ep}$, $D_{\eta_\ep}$, $\nabla_{\eta_\ep}$ and ${\rm div}_{\eta_\ep}$ by using rescaling (\ref{dilatacion_t}) in the set $I_\ep$. The definitions are analogous to the operators  $\mathbb{D}_{h_\ep}$, $\Delta_{h_\ep}$, $D_{h_\ep}$, $\nabla_{h_\ep}$ and ${\rm div}_{h_\ep}$,  just replacing $h_\ep$ by $\eta_\ep$.

 Let $C^\infty_{\rm per}(Y)$ be the space of infinitely differentiable functions in $\mathbb{R}^2$ that are $Y$-periodic. By $L^r_{\rm per}(Y)$ (resp. $W^{1,r}_{\rm per}(Y)$) we denote its completion in the norm $L^r(Y)$ (resp. $W^{1,r}(Y)$). We denote by $L^{r'}_0$ the space of functions of $L^{r'}$ with null integral and by $L^{r'}_{0,{\rm per}}(Y)$  the space of functions in $L^{r'}_{\rm per}(Y)$ with zero mean value.

We denote by $:$ the full contraction of two matrices, i.e. for $A=(a_{ij})_{1\leq i,j\leq 3}$ and $B=(a_{ij})_{1\leq i,j\leq 2}$, we have $A:B=\sum_{i,j=1}^2a_{ij}b_{ij}$.  The canonical basis in $\mathbb{R}^2$ is denoted by $\{e_1,e_2\}$.

Finally, we denote by $O_\ep$ a generic real sequence, which tends to zero with $\ep$ and can change from line to line, and by $C$ a generic positive constant which also can change from line to line.

\subsection{Statement of the problem}  Let us consider a sequence $(u_\ep, p_\ep)\in W^{1,r}_0(D_\ep)^2\times L^{r'}_0(D_\ep)$, $1<r<+\infty$, which satisfies 
\begin{equation}
\left\{
\begin{array}
[c]{r@{\;}c@{\;}ll}%
\\
\displaystyle -\nu\,{\rm div}\left(|\mathbb{D}[u_{\varepsilon}]|^{r-2}\mathbb{D}[u_{\varepsilon }] \right)+ \nabla p_{\varepsilon } &
= &
f\quad \hbox{in }D_\varepsilon,\\
\\
{\rm div}(u_{\varepsilon}) & = & 0 \quad \hbox{in }D_\varepsilon,
\end{array}
\right. \label{system_1_dimension_1}%
\end{equation}  
and the boundary condition
\begin{equation}\label{system_1_dimension_1bc} u_\varepsilon=0\quad \hbox{on }\partial T_\varepsilon \cup\partial \Lambda_\ep,
\end{equation}
where $\nu>0$ is the consistency, and $r'=r/(r-1)$ is the conjugate exponent of $r$. We assume 
\begin{equation}\label{fep}
 f=(f_1(x_1),0)\quad\hbox{with}\quad f_1\in L^{r'}(\omega),
\end{equation}
  which is usual when we deal with thin domains. Since the thickness of the domain  is small, then the vertical component of the force can be neglected, and moreover, the force can be considered independent of the vertical variable.
\\

Our aim is to describe the asymptotic behavior of the velocity $u_\ep$ and the pressure
$p_\ep$ of the fluid as $\ep$ tends to zero and identify an homogenized model coupling the effects of the thickness and the microgeometry of the domain.  To do this,  we will use the equivalent weak variational formulation of (\ref{system_1_dimension_1})--(\ref{system_1_dimension_1bc}), which is the following one: find $u_\ep \in W^{1,r}_0(D_\ep)^2$ and $p_\ep\in L^{r'}_0(D_\ep)$ such that
\begin{equation}\label{form_var_1}
 \begin{array}{l}
\displaystyle \nu\int_{D_\ep}|\mathbb{D}[u_\ep]|^{r-2}\mathbb{D}[u_\ep]:\mathbb{D}[\varphi]\,dx-\int_{D_\ep}p_\ep\,{\rm div}(\varphi)\,dx=\int_{D_\ep}f\cdot\varphi\,dx,
\quad\forall\,\varphi  \in   W^{1,r}_0(D_\ep)^2,\\
\\
\displaystyle 
\int_{D_\ep}{\rm div}(u_\ep)\,\psi\,dx=0\quad\forall\,\psi  \in  L^{r'}(D_\ep).
\end{array} 
\end{equation}
It is well known (see for instance the classical theory \cite{Temam}) that, for every $\ep>0$,  problem (\ref{system_1_dimension_1})--(\ref{system_1_dimension_1bc}) has a unique weak solution $(u_\ep, p_\ep)\in W^{1,r}_0(D_\ep)^2\times L^{r'}_0(D_\ep)$.\\
 
In order to find the limit problem when $\ep$ tends to zero, it is necessary to obtain a priori estimates in fixed domains (with respect to $\ep$), so we introduce the rescaling given by (\ref{dilatacion_p}) for the thin porous media and (\ref{dilatacion_t}) for the thin film, that is
\begin{equation}\label{dilatation_both}\left\{\begin{array}{l}
\displaystyle z_1=x_1, \quad z_2={x_2\over h_\ep}\quad \hbox{if } (x_1, x_2)\in \Omega_\ep,\\
\\
\displaystyle  z_1=x_1, \quad z_2={x_2\over \eta_\ep}\quad \hbox{if } (x_1, x_2)\in I_\ep.
\end{array}\right.
\end{equation}
Using this rescaling, we can define  $\widetilde u_\ep \in W^{1,r}_0(\widetilde D_\ep)^2$ and $\widetilde p_\ep\in L^{r'}_0(\widetilde D_\ep)$ by 
\begin{equation}\label{unk_dilat}
\left\{\begin{array}{l}
\displaystyle \widetilde u_\ep(z)=u_\ep(z_1, h_\ep z_2),\quad  \widetilde p_\ep(z)=p_\ep(z_1, h_\ep z_2)\quad\hbox{if }z\in \widetilde\Omega_\ep\,,
\\
\\
\displaystyle \widetilde u_\ep(z)=u_\ep(z_1, \eta_\ep z_2),\quad  \widetilde p_\ep(z)=p_\ep(z_1, \eta_\ep z_2)\quad\hbox{if }z\in \widetilde I_1\,,
\end{array}\right.
\end{equation}
so the rescaled weak variational formulation is the following: find $\widetilde u_\ep \in W^{1,r}(\widetilde D_\ep)^2$, $\widetilde p_\ep\in L^{r'}_0(\widetilde D_\ep)$ such that
\begin{equation}\label{form_var_1_tilde}
 \begin{array}{l}
\displaystyle \nu\int_{\widetilde \Omega_\ep}h_\ep|\mathbb{D}_{h_\ep}[\widetilde u_\ep]|^{r-2}\mathbb{D}_{h_\ep}[\widetilde u_\ep]:\mathbb{D}_{h_\ep}[\widetilde \varphi]\,dz+\nu\int_{\widetilde I_1}\eta_\ep|\mathbb{D}_{\eta_\ep}[\widetilde u_\ep]|^{r-2}\mathbb{D}_{\eta_\ep}[\widetilde u_\ep]:\mathbb{D}_{\eta_\ep}[\widetilde \varphi]\,dz\\
\\
\displaystyle  -\int_{\widetilde \Omega_\ep}h_\ep\widetilde p_\ep\,{\rm div}_{h_\ep}(\widetilde \varphi)\,dz-\int_{\widetilde I_1}\eta_\ep \widetilde p_\ep\,{\rm div}_{\eta_\ep}(\widetilde \varphi)\,dz\\
\\
\displaystyle =\int_{\widetilde \Omega_\ep}h_\ep f\cdot\widetilde \varphi\,dz+\int_{\widetilde I_1}\eta_\ep f\cdot\widetilde \varphi\,dz,
\qquad\forall\,\widetilde \varphi  \in   W^{1,r}_0(\widetilde D_\ep)^2,\\
\\
\displaystyle 
\int_{\widetilde \Omega_\ep}h_\ep{\rm div}_{h_\ep}(\widetilde u_\ep)\,\widetilde \psi\,dz+\int_{\widetilde I_1}\eta_\ep{\rm div}_{\eta_\ep}(\widetilde u_\ep)\,\widetilde \psi\,dz=0,\qquad\forall\,\widetilde \psi  \in  L^{r'}(\widetilde D_\ep).
\end{array} 
\end{equation}
Moreover, to study the behavior of the velocity in the two media, we introduce the following notation for velocity: 
$$u_\ep=v_\ep+\mathcal{U}_\ep,$$
where
\begin{itemize}
\item[--] $v_\ep$ denotes the velocity in the thin porous medium, extended by zero to the thin domain $I_\ep$, i.e.
\begin{equation}\label{Funcvep}
\displaystyle v_\ep(x)=\left\{\begin{array}{cl}
u_\ep(x) &\hbox{if }x\in \Omega_\ep,\\
\\
0 &\hbox{if }x\in I_\ep.
\end{array}\right.
\end{equation}
We denote by $\widetilde v_\ep$ the dilated velocity in $\widetilde\Omega_\ep$ obtained from $v_\ep$ by using the change of variables (\ref{dilatacion_p}).  Then,  $\widetilde v_\ep$ satisfies the following equality
\begin{equation}\label{form_var_1_tildev}
 \begin{array}{l}
\!\!\!\!\!\!\displaystyle \nu\int_{\widetilde \Omega_\ep}|\mathbb{D}_{h_\ep}[\widetilde v_\ep]|^{r-2}\mathbb{D}_{h_\ep}[\widetilde v_\ep]:\mathbb{D}_{h_\ep}[\widetilde \varphi]\,dz -\int_{\widetilde \Omega_\ep}\widetilde p_\ep\,{\rm div}_{h_\ep}(\widetilde \varphi)\,dz=\int_{\widetilde \Omega_\ep}f\cdot\widetilde \varphi\,dz,
\quad\forall\,\widetilde \varphi  \in   W^{1,r}_0(\widetilde \Omega_\ep)^2,\\
\\
\displaystyle 
\!\!\!\!\!\!\int_{\widetilde \Omega_\ep}{\rm div}_{h_\ep}(\widetilde v_\ep)\,\widetilde \psi\,dz=0,\quad\forall\,\widetilde \psi  \in  L^{r'}(\widetilde \Omega_\ep).
\end{array} 
\end{equation}

\item[--] $\mathcal{U}_\ep$ denotes the velocity in the thin film $I_\ep$, extended by zero to the thin porous medium $\Omega_\ep$, i.e.
\begin{equation}\label{Funcvep}
\displaystyle \mathcal{U}_\ep(x)=\left\{\begin{array}{cl}
0 &\hbox{if }x\in  \Omega_\ep\,,\\
\\
u_\ep(x) &\hbox{if }x\in  I_\ep\,.
\end{array}\right.
\end{equation}
We denote by $\mathcal{\widetilde U}_\ep$ the dilated velocity in $\widetilde I_1$ obtained from   $\mathcal{U}_\ep$ by using the change of variables (\ref{dilatacion_t}).   Then $\mathcal{\widetilde U}_\ep$ satisfies the following equality
\begin{equation}\label{form_var_1_tildeU}
 \begin{array}{l}
\displaystyle \nu\int_{\widetilde I_1}|\mathbb{D}_{\eta_\ep}[\mathcal{\widetilde U}_\ep]|^{r-2}\mathbb{D}_{\eta_\ep}[\mathcal{\widetilde U}_\ep]:\mathbb{D}_{\eta_\ep}[\widetilde \varphi]\,dz-\int_{\widetilde I_1}\widetilde p_\ep\,{\rm div}_{\eta_\ep}(\widetilde \varphi)\,dz=\int_{\widetilde I_1}f\cdot\widetilde \varphi\,dz,
\quad\forall\,\widetilde \varphi  \in   W^{1,r}_0(\widetilde I_1)^2,\\
\\
\displaystyle 
\int_{\widetilde I_1}{\rm div}_{\eta_\ep}(\mathcal{\widetilde U}_\ep)\,\widetilde \psi\,dz=0\quad\forall\,\widetilde \psi  \in  L^{r'}(\widetilde I_1).
\end{array} 
\end{equation}
\end{itemize}
\section{A priori estimates}\label{sec:estimates}
\subsection{Some technical estimates}
Let us begin with the  Poincar\'e and Korn inequalities in the thin porous medium $\Omega_\ep$.
\begin{lemma}[Poincare and Korn inequalities] \label{Lemma_Poincare}  For every $\varphi\in W^{1,r}(\Omega_{\varepsilon})^2$,  $1<r<+\infty$, wtih $\varphi=0$ on $\partial\Omega_\varepsilon \setminus \Sigma$, there exists a positive constant $C$, independent of $\ep$, such that, 
\begin{equation}\label{Poincare1}
\|\varphi\|_{L^r(\Omega_\ep)^2}\leq C\ep\|D \varphi\|_{L^r(\Omega_\ep)^{2\times 2}},
\end{equation}
\begin{equation}\label{Korn1}
\left\Vert D \varphi\right\Vert_{L^r(\Omega_{\varepsilon})^{2 \times 2}}\leq C\left\Vert \mathbb{D}[\varphi]\right\Vert_{L^r(\Omega_{\varepsilon})^{2\times 2}}.\end{equation}
\end{lemma}
 
\begin{proof} We observe that $\Omega_\ep$ can be divided in small cubes of lateral and vertical length $\ep$. We consider the periodic cell $Y_f$  and  have a Friedrichs inequality    
\begin{equation}\label{proof_gaffney}
\int_{Y_f}|\varphi|^r\,dz\leq C\int_{Y_f}|D\varphi|^r\,dz\,,
\end{equation}
for every $\varphi(z)\in W^{1,r}(Y_f)^2$ such that $\varphi =0$ on $\partial T$, where the constant $C$ depends only on $Y_f$. Then, for every $k\in \mathbb{Z}^2$, by the change of variable
\begin{equation}\label{change}
\begin{array}{l}
\displaystyle k+z={x\over \ep},\quad dz={dx\over \ep^2},\quad  \partial_{z}=\ep\partial_{x},\end{array}
\end{equation}
we rescale (\ref{proof_gaffney}) from $Y_f$ to $Y_{f_k,\ep}$.  This yields that, for every function $\varphi(x)\in W^{1,r}(Y_{f_k,\ep})^2$, one has
$$
\int_{Y_{f_k,\ep}}|\varphi|^r\,dx\leq C\varepsilon^r\int_{Y_{f_k,\ep}}|D\varphi|^r\,dx\,,
$$
with the same constant $C$ as in (\ref{proof_gaffney}).  Summing previous inequality for every $k\in \mathcal{K}_\ep$,  we get (\ref{Poincare1}). 

Finally, Korn's inequality (\ref{Korn1}) follows from the classical Korn inequality, see \cite{Boyer}.

\end{proof}

Next, we give an useful estimate in the thin film $I_\ep$.

\begin{lemma}\label{Lem:Iep}
For every function $\varphi\in W^{1,r}_0(D_\ep)^2$, with $1<r<\infty$, there exists a constant $C>0$ independent of $\ep$, such that,
\begin{equation}\label{estim_fissure_1}
\|\varphi\|_{L^r(I_\ep)^2}\leq C\eta_\ep^{1\over 2}(\eta_\ep+\ep)^{1\over 2}\|\mathbb{D}[\varphi]\|_{L^r(D_\ep)^{2\times 2}}.
\end{equation}
\end{lemma}

\begin{proof} Because the thickness of  $I_\ep$ (see for instance \cite{MikelicTapiero}), we have 
\begin{equation}\label{PoincareIep}\|\varphi\|_{L^r(I_\ep)^2}\leq C\eta_\ep\|D\varphi\|_{L^r(I_\ep)^{2\times 2}}.
\end{equation}
Next, if we choose a point $t\in T_\ep$, which is close to the point $x\in I_\ep$, then, we have 
$$|\varphi(x)-\varphi(t)|=|D\varphi(\xi) (x-t)|\leq (\ep+\eta_\ep)|D\varphi|.$$
Since $\varphi(t)=0$ because $t\in T_\ep$, we have
$$\|\varphi\|_{L^r(I_\ep)^2}\leq C(\ep+\eta_\ep)\|D\varphi\|_{L^r(I_\ep)^{2\times 2}}.$$
Multiplying the above inequality with inequality (\ref{PoincareIep}), we get
\begin{equation}\label{PoincareIep1}
\|\varphi\|_{L^r(I_\ep)^2}\leq C\eta_\ep^{1\over 2}(\ep+\eta_\ep)^{1\over 2}\|Dv\|_{L^r(I_\ep)^{2\times 2}}\leq C\eta_\ep^{1\over 2}(\ep+\eta_\ep)^{1\over 2}\|Dv\|_{L^r(D_\ep)^{2\times 2}},
\end{equation}
and from the classical Korn inequality in $I_\ep$, we obtain the estimate (\ref{estim_fissure_1}).
\end{proof}

\subsection{Estimates for velocity} 
We derive the estimates for velocity in the whole domain $D_\ep$ for $u_\ep$, and also, in the sets $\Omega_\ep$ and $I_\ep$ for $v_\ep$ and $\mathcal{U}_\ep$ respectively.

\begin{lemma}\label{Lem:estimates_velocity}There exists a constant $C>0$ independent of $\ep$, such that if $u_\ep\in W^{1,r}_0(D_\ep)^2$, with $1<r<+\infty$, is the solution of problem (\ref{system_1_dimension_1})--(\ref{system_1_dimension_1bc}), it holds
\begin{equation}\label{estim_u_O}
\|v_\ep\|_{L^r(\Omega_\ep)^{2}}\leq C\left( \eta_\ep^{2r-1\over r}\ep^{r-1}+h_\ep^{r-1\over r}\ep^r\right)^{1\over r-1},
\end{equation}
\begin{equation}\label{estim_u_I}
\|\mathcal{U}_\ep\|_{L^r(I_\ep)^{2}}\leq C
\eta_\ep^{1+{2r-1\over r(r-1)}}+\ep^{1\over r-1} \eta_\ep  h_\ep^{1\over r}
+ \eta_\ep^{1\over 2}\ep^{{1\over 2}+{1\over r-1}}h_\ep^{1\over r},
\end{equation}
\begin{equation}\label{estim_DDu}
\|\mathbb{D}[u_\ep]\|_{L^r(D_\ep)^{2\times 2}}\leq C\left( \eta_\ep^{2r-1\over r}+h_\ep^{r-1\over r}\ep\right)^{1\over r-1},
\end{equation}
\begin{equation}\label{estim_Du}
\|Du_\ep\|_{L^r(D_\ep)^{2\times 2}}\leq C\left( \eta_\ep^{2r-1\over r}+h_\ep^{r-1\over r}\ep\right)^{1\over r-1}.
\end{equation}
\end{lemma}

\begin{proof}
Using $u_\ep$ as test function in (\ref{form_var_1}), we have
\begin{equation}\label{estim_proof_1}
\nu\|\mathbb{D}[u_\ep]\|_{L^r(D_\ep)^{2\times 2}}^r=\int_{D_\ep}f\cdot u_\ep\,dx.
\end{equation}
Using the H${\rm \ddot{o}}$lder inequality and the assumption of $f$ given in (\ref{fep}), we obtain that there exists a constant $C>0$ such that
$$\int_{D_\ep}f\cdot u_\ep\,dx\leq C\left(\eta_\ep^{1\over r'}\|u_\ep\|_{L^r(I_\ep)^2}+h_\ep^{1\over r'}\|u_\ep\|_{L^r(\Omega_\ep)^2}\right),$$
and by inequalities (\ref{Poincare1}), (\ref{Korn1}) and (\ref{estim_fissure_1}), we have
$$\begin{array}{ll}
\displaystyle \int_{D_\ep}f\cdot u_\ep\,dx &\displaystyle \leq C\left( \eta_\ep^{1\over r'}\eta_\ep^{1\over 2}(\ep+\eta_\ep)^{1\over 2}+h_\ep^{1\over r'}\ep\right)\|\mathbb{D}[u_\ep]\|_{L^r(D_\ep)^{2\times 2}}
\\
\\
&\displaystyle \leq C\left( \eta_\ep^{{1\over r'}+1}+\eta_\ep^{{1\over r'}+{1\over 2}}\ep^{1\over 2}+h_\ep^{1\over r'}\ep\right)\|\mathbb{D}[u_\ep]\|_{L^r(D_\ep)^{2\times 2}}.
\end{array}$$
Therefore, from equation (\ref{estim_proof_1}), we get
$$\|\mathbb{D}[u_\ep]\|_{L^r(D_\ep)^{2\times 2}}\leq C\left( \eta_\ep^{{1\over r'}+1}+\eta_\ep^{{1\over r'}+{1\over 2}}\ep^{1\over 2}+h_\ep^{1\over r'}\ep\right)^{1\over r-1}.$$
Since $\ep\ll\eta_\ep$, then $\eta_\ep^{{1\over r'}+{1\over 2}}\ep^{1\over 2}\ll \eta_\ep^{{1\over r'}+1}$ and so, the term $\eta_\ep^{{1\over r'}+{1\over 2}}\ep^{1\over 2}$ can be dropped. Then, taking into account that $1/r'+1=(2r-1)/r$ and $1/r'=(r-1)/r$, we get estimate (\ref{estim_DDu}). From the  classical Korn inequality, we have inequality (\ref{estim_Du}). Applying inequality (\ref{Poincare1}) together with inequality (\ref{estim_Du}), we obtain inequality (\ref{estim_u_O}).

\noindent Finally, applying inequalities (\ref{estim_fissure_1}) and (\ref{estim_DDu}), we get
$$\begin{array}{rl}
\|\mathcal{U}_\ep\|_{L^r(I_\ep)^2}\leq &\displaystyle C(\eta_\ep+\eta_\ep^{1\over 2}\ep^{1\over 2})\left( \eta_\ep^{2r-1\over r(r-1)}+h_\ep^{1\over r}\ep^{1\over r-1}\right)\\
\\
= &\displaystyle C\left(\eta_\ep^{1+{2r-1\over r(r-1)}}+\ep^{1\over r-1}\eta_\ep h_\ep^{1\over r}+\eta_\ep^{1\over 2}\eta_\ep^{2r-1\over r(r-1)}\ep^{1\over 2}+\eta_\ep^{1\over 2}\ep^{{1\over 2}+{1\over r-1}}h_\ep^{1\over r}\right).
\end{array}$$
Since $\ep\ll \eta_\ep$, then $\eta_\ep^{1\over 2}\eta_\ep^{2r-1\over r(r-1)}\ep^{1\over 2}\ll\eta_\ep^{1+{2r-1\over r(r-1)}}$   and so, the term $\eta_\ep^{1\over 2}\eta_\ep^{2r-1\over r(r-1)}\ep^{1\over 2}$ can be dropped, and inequality (\ref{estim_u_I}) holds.
\end{proof}
As consequence, by using the change of variables (\ref{dilatacion_p}) in the thin porous medium $\Omega_\ep$ and (\ref{dilatacion_t}) in the thin film $I_\ep$, we derive the estimates for the rescaled velocities.
\begin{corollary}\label{estimates_vel_tilde}There exists a constant $C>0$ independent of $\ep$, such that we have the following estimates depending on the media:
\begin{itemize}
\item[--] In the porous media $\widetilde \Omega_\ep$, we have
\begin{equation}\label{estim_u_O_dil}
\|\widetilde v_\ep\|_{L^r(\widetilde \Omega_\ep)^{2}}\leq C\left( \eta_\ep^{2r-1\over r}\ep^{r-1}h_\ep^{-{r-1\over r}}+ \ep^r\right)^{1\over r-1},
\end{equation}
\begin{equation}\label{estim_DDu_O_tilde}
\|\mathbb{D}_{h_\ep}[\widetilde v_\ep]\|_{L^r(\widetilde \Omega_\ep)^{2\times 2}}\leq C\left( \eta_\ep^{2r-1\over r}h_\ep^{-{r-1\over r}}+\ep\right)^{1\over r-1},
\end{equation}
\begin{equation}\label{estim_Du_O_dil}
\|D_{h_\ep}\widetilde v_\ep\|_{L^r(\widetilde \Omega_\ep)^{2\times 2}}\leq C\left( \eta_\ep^{2r-1\over r}h_\ep^{-{r-1\over r}}+\ep\right)^{1\over r-1}.
\end{equation}
\item[--] In the  the free media $\widetilde I_1$, we have
\begin{equation}\label{estim_u_I_dil}
\|\mathcal{\widetilde U}_\ep\|_{L^r(\widetilde I_1)^{2}}\leq C \left( \eta_\ep^{r\over r-1}+\ep^{1\over r-1}h_\ep^{1\over r}\eta_\ep^{r-1\over r}+\ep^{r+1\over 2(r-1)}h_\ep^{1\over r}\eta_\ep^{r-2\over 2r}\right),
\end{equation}
\begin{equation}\label{estim_DDu_I_dil}
\|\mathbb{D}_{\eta_\ep}[\mathcal{\widetilde U}_\ep]\|_{L^r(\widetilde I_1)^{2\times 2}}\leq C\left( \eta_\ep+h_\ep^{r-1\over r}\ep\eta_\ep^{-{r-1\over r}}\right)^{1\over r-1},
\end{equation}
\begin{equation}\label{estim_Du_I_dil}
\|D_{\eta_\ep}\mathcal{\widetilde U}_\ep\|_{L^r(\widetilde I_1)^{2\times 2}}\leq C\left( \eta_\ep+h_\ep^{r-1\over r}\ep\eta_\ep^{-{r-1\over r}}\right)^{1\over r-1}.
\end{equation}

\end{itemize}
\end{corollary}
\begin{proof} Estimates for dilated velocity (\ref{estim_u_O_dil})--(\ref{estim_Du_O_dil})  in $\widetilde \Omega_\ep$ are obtained directly from (\ref{estim_u_O}), (\ref{estim_DDu}) and (\ref{estim_Du}) by applying the change of variables (\ref{dilatacion_p}), just taking into account that
$$\begin{array}{c}
\displaystyle \|v_\ep\|_{L^r(\Omega_\ep)^2}=h_\ep^{{1\over r}}\|\widetilde v_\ep\|_{L^r(\widetilde \Omega_\ep)^2},\  \|\mathbb{D} v_\ep\|_{L^r(\Omega_\ep)^{2\times 2}}=h_\ep^{{1\over r}}\|\mathbb{D}_{h_\ep} \widetilde v_\ep\|_{L^r(\widetilde \Omega_\ep)^{2\times 2}},\\
\displaystyle   \|D  v_\ep\|_{L^r(  \Omega_\ep)^{2\times 2}}=h_\ep^{{1\over r}}\|D_{h_\ep} \widetilde v_\ep\|_{L^r(\widetilde \Omega_\ep)^{2\times 2}}.
\end{array}$$
Similarly, estimates for dilated velocity (\ref{estim_u_I_dil})--(\ref{estim_Du_I_dil})  in $\widetilde I_1$ are obtained directly from (\ref{estim_u_I}), (\ref{estim_DDu}) and (\ref{estim_Du}) by applying the change of variables (\ref{dilatacion_t}), just taking into account that
$$\begin{array}{c}\|\mathcal{U}_\ep\|_{L^r(I_\ep)^2}=\eta_\ep^{{1\over r}}\|\mathcal{\widetilde U}_\ep\|_{L^r(\widetilde I_1)^2},\  \|\mathbb{D} [\mathcal{U}_\ep]\|_{L^r(I_\ep)^{2\times 2}}=\eta_\ep^{{1\over r}}\|\mathbb{D}_{\eta_\ep}[\mathcal{\widetilde U}_\ep]\|_{L^r(\widetilde I_1)^{2\times 2}},\\
\displaystyle  \|D  \mathcal{U}_\ep\|_{L^r(  I_\ep)^{2\times 2}}=\eta_\ep^{{1\over r}}\|D_{\eta_\ep} \mathcal{\widetilde U}_\ep\|_{L^r(\widetilde I_1)^{2\times 2}}.
\end{array}$$
\end{proof}

\subsection{Estimates for pressure in the porous part}\label{sec:restriction}
 Next, we derive a priori estimates for the pressure in the porous part. To do this, we need to extend the pressure to the whole thin film $Q_\ep$ (which also depends on $\ep$).  To do this, we  generalize a result from  \cite[Lemma 3.3]{BayadaThinThin}  (see also \cite[Lemma 5.3]{SG_MANA}) by introducing a restriction operator $\mathcal{R}^\ep_r$ from $W^{1,r}_0(Q_\ep)^2$ into $W^{1,r}_0(\Omega_\ep)^2$, $1<r<+\infty$. We remark that in the case $r=2$, this restriction operator $\mathcal{R}^\ep_2$ agrees with the one defined in \cite[Lemma 3.3]{BayadaThinThin} and in \cite[Lemma 5.3]{SG_MANA}.
\begin{lemma} \label{restriction_operator}
There exists  a (restriction) operator $\mathcal{R}^\ep_r$ acting from $W^{1,r}_0(Q_\ep)^2$ into $W^{1,r}_0(\Omega_\ep)^2$ such that
\begin{enumerate}
\item $\mathcal{R}^\ep_r \varphi=\varphi$, if $\varphi \in W^{1,r}_0(\Omega_\ep)^2$.
\item ${\rm div}(\mathcal{R}^\ep_r \varphi)=0\hbox{  in }\Omega_\ep$, if ${\rm div}(\varphi)=0\hbox{  on }Q_\ep$.
\item For every $\varphi\in W^{1,r}_0(Q_\ep)^3$, there exists a positive constant $C$, independent of $\varphi$ and $\ep$, such that
\begin{equation}\label{estim_restricted}
\begin{array}{l}
\|\mathcal{R}^\ep_r \varphi\|_{L^r(\Omega_\ep)^{2}}+ \ep\|D \mathcal{R}^\ep_r \varphi\|_{L^r(\Omega_\ep)^{2\times 2}} \leq C\left(\|\varphi\|_{L^r(Q_\ep)^2}+\ep \|D\varphi\|_{L^r(Q_\ep)^{2\times 2}}\right)\,.
\end{array}
\end{equation}
\end{enumerate}
\end{lemma}
\begin{proof}
  Let us consider the linear map $\mathcal{R}_r$ constructed in \cite[Lemma 1.1]{Bourgeat1} from $W^{1,r}(Y)^2\to W^{1,r}_T(Y_f)^2$, $1<r<+\infty$, where $W^{1,r}_T(Y_f)^2=\{\varphi\in W^{1,r}(Y_f)^2\ :\ \varphi=0\hbox{ on }T\}$, such that
\begin{equation}\label{Rr}
\|\mathcal{R}_r\varphi\|_{W^{1,r}(Y_f)^2}\leq C\|\varphi\|_{W^{1,r}(Y)^2},
\end{equation}
and $\mathcal{R}_r\varphi$ coincides with $\varphi$ if $\varphi$ is zero on $T$ (i.e. if $\varphi\in W^{1,r}_T(Y_f)^2$)  and ${\rm div}(\varphi)=0$ implies ${\rm div}(\mathcal{R}_r\varphi)=0$. Then, $\mathcal{R}^\ep_r$ is defined by applying $\mathcal{R}_r$ to each $Y_{k,\ep}$. Consequently, the two first items are satisfied. 

Finally, we will prove the third item. From (\ref{Rr}), by the change of variables (\ref{change}), as in Lemma \ref{Lemma_Poincare}, we have
$$\begin{array}{l}
\displaystyle
\int_{Y_{f_k,\varepsilon}}|\mathcal{R}^\ep_r\varphi|^r\,dx+\ep^r\int_{Y_{f_k,\varepsilon}}|D\mathcal{R}_r^\ep \varphi|^r\,dx \leq C\left(\int_{Y_{k,\varepsilon}}|\varphi|^r\,dx+\ep^r\int_{Y_{k,\varepsilon}}|D\varphi|^r\,dx
\right).
\end{array}
$$
and so, summing previous inequality for every $k\in \mathcal{K}_\ep$, we deduce (\ref{estim_restricted}).
\end{proof}

\begin{lemma}\label{restriction_operator2} Setting $\mathcal{\widetilde R}_r^\ep\widetilde\varphi=\mathcal{R}_r^\ep\varphi$ for any $\widetilde \varphi\in W^{1,r}_0(\Omega)^2$, $1<r<+\infty$, where $\widetilde \varphi$ is obtained from $ \varphi$ by using the change of variables (\ref{dilatacion_p}),  and $\mathcal{R}_r^\ep$ is defined in Lemma \ref{restriction_operator}, we have the following estimates:
\begin{equation}\label{estim_restricted2}
\begin{array}{l}
\displaystyle 
\|\mathcal{\widetilde R}^\ep_r \widetilde \varphi\|_{L^r(\widetilde \Omega_\ep)^{2}}\leq C\|\widetilde \varphi\|_{W^{1,r}_0(\Omega)^2},\quad \|D_{h_\ep} \mathcal{\widetilde R}^\ep_r \widetilde \varphi\|_{L^r(\widetilde \Omega_\ep)^{2\times 2}} \leq C\ep^{-1}\|\widetilde \varphi\|_{W^{1,r}_0(\Omega)^2}.
\end{array}
\end{equation}

\end{lemma}
\begin{proof}
Applying the change of variables (\ref{dilatacion_p}) to estimates (\ref{estim_restricted}) and taking into account that $\ep\ll h_\ep$, we get 
$$\begin{array}{rl}
\displaystyle
\|\mathcal{\widetilde R}^\ep_r \widetilde \varphi\|_{L^r(\widetilde \Omega_\ep)^{2}}+\ep\|D_{h_\ep}\mathcal{\widetilde R}^\ep_r \widetilde \varphi\|_{L^r(\widetilde \Omega_\ep)^{2\times 2}} \leq & \displaystyle C\left(\|\widetilde \varphi\|_{L^r(\Omega)^2}+\ep \|D_{h_\ep}\widetilde \varphi\|_{L^r(\Omega)^{2\times 2}}\right)\\
\\ \leq &
\displaystyle\displaystyle C\left( \|D\widetilde \varphi\|_{L^r(\Omega)^{2\times 2}}+\ep h_\ep^{-1}\|D\widetilde \varphi\|_{L^r(\Omega)^{2\times 2}}\right)\\
\\
\leq &\displaystyle
C \|\widetilde \varphi\|_{W^{1,r}_0(\Omega)^2},
\end{array}$$
which implies estimates (\ref{estim_restricted2}).
\end{proof}

Denoting by $p^1_\ep$ the restriction to $\Omega_\ep$ of the overall pressure $p_\ep$,  with  the additive constant being determined by $\int_{\Omega_\ep}p^1_\ep\,dx=0$, we give the existence of an extended pressure to $Q_\ep$ by duality arguments. 
\begin{lemma} \label{Estimates_extended_lemma}  There exists an extension $P^1_\ep\in L^{r'}_0(Q_\ep)$ of the pressure $p_\ep^1$. Moreover, defining the dilated and extended pressure $\widetilde P^1_\ep\in L^{r'}_0(\Omega)$ obtained from $P^1_\ep$ by using the change of variables (\ref{dilatacion_p}), then 
there exists a positive constant $C$ independent of $\ep$, such that 
\begin{equation}\label{esti_P}
\|\widetilde P_\ep^1\|_{L^{r'}(\Omega)}\leq C\left( h_\ep^{-{r-1\over r}}\eta_\ep^{2r-1\over r}\ep^{-1}+1\right),\quad 
\|\nabla_{h_\ep}  \widetilde P_\ep^1\|_{W^{-1,r'}(\Omega)^2}\leq C\left( h_\ep^{-{r-1\over r}}\eta_\ep^{2r-1\over r}\ep^{-1}+1\right).
\end{equation}
\end{lemma}
\begin{proof}  We divide the proof in two steps.

{\it Step 1. Extension of $p_\ep^1$ to $Q_\ep$}. Using the restriction operator $\mathcal{R}^\ep_r$ given in  Lemma \ref{restriction_operator}, we   introduce $F_\ep$ in $W^{-1,r'}(Q_\ep)^2$ in the following way
\begin{equation}\label{F}\langle F_\varepsilon, \varphi\rangle_{W^{-1,r'}(Q_\varepsilon)^2, W^{1,r}_0(Q_\ep)^2}=\langle \nabla p_\varepsilon, \mathcal{R}^\varepsilon_r \varphi\rangle_{{W^{-1,r'}(\Omega_\varepsilon)^2, W^{1,r}_0(\Omega_\ep)^2}}\,,\quad \hbox{for any }\varphi\in W^{1,r}_0(Q_\varepsilon)^2\,,
\end{equation}
and compute the right hand side of (\ref{F}) by using  in (\ref{form_var_1}), which  gives
\begin{equation}\label{equality_duality}
\begin{array}{l}
\displaystyle
\left\langle F_{\varepsilon},\varphi\right\rangle_{W^{-1,r'}(Q_\varepsilon)^2, W^{1,r}_0(Q_\ep)^2}=\displaystyle
-\nu\int_{\Omega_\varepsilon}|\mathbb{D}[v_\ep]|^{r-2}\mathbb{D}[v_\ep] : \mathbb{D}[\mathcal{R}^{\varepsilon}_r\varphi]\,dx+ \int_{\Omega_\varepsilon} f\cdot (\mathcal{R}^{\varepsilon}_r\varphi)\,dx \,.
\end{array}\end{equation}
Using Corollary \ref{estimates_vel_tilde} for fixed $\ep$, we see that it is a bounded functional on $W^{1,r}_0(Q_\ep)$ (see Step  2 of the proof), and in fact $F_\ep\in W^{-1,r'}(Q_\ep)^3$. Moreover, ${\rm div}(\varphi)=0$ implies $\left\langle F_{\varepsilon},\varphi\right\rangle=0\,,$ and the DeRham theorem gives the existence of $P^1_\varepsilon\in L^{r'}_0(Q_\varepsilon)$ with $F_\varepsilon=\nabla P^1_\varepsilon$.
\\

{\it Step 2. Estimates for dilated and extended pressure.} Consider $\widetilde P^1_\ep$ obtained from $P^1_\ep$ by using the change of variables (\ref{dilatacion_p}). By using the Ne${\breve{\rm c}}$as inequality (see for instance \cite{Boyer}) for $\widetilde P^1_\ep\in L^{r'}_0(\Omega)$, then 
$$\|\widetilde P^1_\ep\|_{L^{r'}(\Omega)}\leq C\|\nabla_z  \widetilde P^1_\ep\|_{W^{-1,r'}(\Omega)^2}\leq C\|\nabla_{h_\ep} \widetilde P^1_\ep\|_{W^{-1,r'}(\Omega)^2},$$
and thus, to prove (\ref{esti_P}), it is enough to prove the second estimate  in (\ref{esti_P}) for $\nabla_{h_\ep} \widetilde P_\ep^1$. \\

Let us prove it. For any $\widetilde\varphi\in W^{1,r}_0(\Omega)^2$, using the change of variables (\ref{dilatacion_p}), we have
$$\begin{array}{rl}\left\langle \nabla_{h_\ep}\widetilde P^1_\ep, \widetilde \varphi\right\rangle_{W^{-1,r'}(\Omega)^2, W^{1,r}_0(\Omega)^2}=&\displaystyle -\int_{\Omega}\widetilde P^1_\ep\,{\rm div}_{h_\ep}(\widetilde\varphi)\,dz\\
=&\displaystyle-h_\ep^{-1}\int_{Q_\ep}P^1_\ep\,{\rm div}(\varphi)\,dx=h_\ep^{-1}\left\langle \nabla  P^1_\ep,  \varphi\right\rangle_{W^{-1,r'}(Q_\ep)^2, W^{1,r}_0(Q_\ep)^2}.
\end{array}$$
Then, using the identification (\ref{equality_duality}) of $F_\ep$, we get
$$
\begin{array}{l}
\displaystyle
\left\langle \nabla_{h_\ep}\widetilde P^1_\ep,\widetilde \varphi\right\rangle_{W^{-1,r'}(\Omega)^2, W^{1,r}_0(\Omega)^2}=\displaystyle h_\ep^{-1}\left(
-\nu\int_{\Omega_\varepsilon}|\mathbb{D}[v_\ep]|^{r-2}\mathbb{D}[v_\ep] : \mathbb{D}[\mathcal{R}^{\varepsilon}_r\varphi]\,dx+ \int_{\Omega_\varepsilon} f\cdot (\mathcal{R}^{\varepsilon}_r\varphi)\,dx\right) \,,
\end{array}
$$
and applying the change of variables (\ref{dilatacion_p}), we get
\begin{equation}\label{equality_duality2}
\begin{array}{l}
\displaystyle
\left\langle \nabla_{h_\ep}\widetilde P^1_\ep,\widetilde \varphi\right\rangle_{W^{-1,r'}(\Omega)^2, W^{1,r}_0(\Omega)^2}=\displaystyle  
-\nu\int_{\widetilde \Omega_\varepsilon}|\mathbb{D}_{h_\ep}[\widetilde v_\ep]|^{r-2}\mathbb{D}_{h_\ep}[\widetilde v_\ep] : \mathbb{D}_{h_\ep}[ \mathcal{\widetilde R}^{\varepsilon}_r\widetilde \varphi]\,dz+ \int_{\widetilde \Omega_\varepsilon} f\cdot (\mathcal{\widetilde R}^{\varepsilon}_r\widetilde \varphi)\,dz  \,.
\end{array}
\end{equation}
Let us now estimate the right-hand side of (\ref{equality_duality2}).
From the H${\rm \ddot{o}}$lder inequality and using estimates for $\widetilde v_\ep$ in (\ref{estim_u_O_dil})--(\ref{estim_Du_O_dil}), assumption of $f$ given in (\ref{fep}) and estimates of the dilated restricted operator (\ref{estim_restricted2}),  we   obtain
$$
\begin{array}{rl}
\displaystyle
\left|\int_{ \widetilde \Omega_\ep}|\mathbb{D}_{h_\ep}[\widetilde v_\ep]|^{r-2}\mathbb{D}_{h_\ep}[\widetilde v_\ep] :\mathbb{D}_{h_\ep}[\mathcal{\widetilde R}^\ep_r\widetilde \varphi]\,dz\right|\leq &\displaystyle C\|\mathbb{D}_{h_\ep}[\widetilde v_\ep]\|_{L^r(\widetilde \Omega_\ep)^{2\times 2}}^{r-1}\|D_{h_\ep} \mathcal{\widetilde R}^\ep_r \widetilde \varphi\|_{L^r(\widetilde \Omega_\ep)^{2\times 2}}\\
 \leq& C\left( h_\ep^{-{r-1\over r}}\eta_\ep^{2r-1\over r}+\ep\right)\|D_{h_\ep}\mathcal{\widetilde R}^\ep_r\widetilde \varphi\|_{L^r(\widetilde \Omega_\ep)^{2\times 2}}\\
 \leq&\displaystyle C\left( h_\ep^{-{r-1\over r}}\eta_\ep^{2r-1\over r}\ep^{-1}+1\right)\|\widetilde \varphi\|_{W^{1,r}_0(\Omega)^2},\\
\\\displaystyle
 \left|\int_{\widetilde \Omega_\ep}f\cdot (\mathcal{\widetilde R}^\ep_r  \widetilde \varphi) \,dz\right|\leq&\displaystyle  C\| \mathcal{\widetilde R}^\ep_r \widetilde\varphi \|_{L^r(\widetilde \Omega_\ep)^2}\leq C\|\widetilde \varphi\|_{W^{1,r}_0(\Omega)^2}\,,
\end{array}
$$
which together with (\ref{equality_duality2}) gives $$\left|\left\langle \nabla_{h_\ep}\widetilde P_\ep,\widetilde \varphi\right\rangle_{W^{-1,r'}(\Omega)^2, W^{1,r}_0(\Omega)^2}\right|\leq C\left( h_\ep^{-{r-1\over r}}\eta_\ep^{2r-1\over r}\ep^{-1}+1\right)\|\widetilde \varphi\|_{W^{1,r}_0(\Omega)^2}.$$
This implies the second estimate given in (\ref{esti_P}), which concludes the proof.
\end{proof}

\subsection{Estimates for pressure in the free part}
Now, we obtain estimates of the pressure in the thin film. 
We denote by $\widetilde P_\ep^2$ the restriction to $\widetilde I_1$ of the overall pressure $\widetilde p_\ep$, i.e. 
\begin{equation}\label{P2}\widetilde P_\ep^2(z)= \widetilde p_\ep(z)-\widetilde c_\ep\quad\hbox{if }z\in \widetilde I_1,
\end{equation}
 with the additive constant $\widetilde c_\ep$ determined by
\begin{equation}\label{cep}\widetilde c_\ep={1\over |\widetilde I_1|}\int_{\widetilde I_1}\widetilde p_\ep\,dz.
\end{equation}

\begin{lemma} \label{Estimates_extended_lemma_thin_film}  
There exists a positive constant $C$ independent of $\ep$, such that 
\begin{equation}\label{esti_P_thin_film}
\|\widetilde P_\ep^2\|_{L^{r'}(\widetilde I_1)}\leq C\left( 1+h_\ep^{r-1\over r}\ep\eta_\ep^{-{2r-1\over r}}
 \right),\quad 
\|\nabla_{\eta_\ep}  \widetilde P_\ep^2\|_{W^{-1,r'}(\widetilde I_1)^2}\leq C\left( 1+h_\ep^{r-1\over r}\ep\eta_\ep^{-{2r-1\over r}}
 \right).
\end{equation}
\end{lemma}
\begin{proof}
Let us first remark that, by using the Ne${\breve{\rm c}}$as inequality (see for instance \cite{Boyer}) for $\widetilde P^2_\ep\in L^{r'}_0(\widetilde I_1)$, then 
$$\|\widetilde P^2_\ep\|_{L^{r'}(\widetilde I_1)}\leq C\|\nabla_z  \widetilde P^2_\ep\|_{W^{-1,r'}(\widetilde I_1)^2}\leq C\|\nabla_{\eta_\ep} \widetilde P^2_\ep\|_{W^{-1,r'}(\widetilde I_1)^2},$$
and thus, to prove (\ref{esti_P_thin_film}), it is enough to prove the second estimate  in (\ref{esti_P_thin_film}) for $\nabla_{\eta_\ep} \widetilde P^2_\ep$. \\

Let us prove it. For any $\widetilde\varphi\in W^{1,r}_0(\widetilde I_1)^2$, using the change of variables (\ref{dilatacion_t}), we have
\begin{equation}\label{equality_duality23}
\left\langle \nabla_{\eta_\ep}\widetilde P^2_\ep, \widetilde \varphi\right\rangle_{W^{-1,r'}(\Omega)^2, W^{1,r}_0(\Omega)^2}=-\nu\int_{\widetilde I_1}|\mathbb{D}_{\eta_\ep}[ \mathcal{\widetilde U}_\ep]|^{r-2}\mathbb{D}_{\eta_\ep}[\mathcal{\widetilde  U}_\ep] : \mathbb{D}_{\eta_\ep}[\widetilde \varphi]\,dz+ \int_{\widetilde I_1} f\cdot \widetilde \varphi\,dz.
\end{equation}
Let us now estimate the right-hand side of this equality:
\begin{itemize}
\item[--] From the H${\rm \ddot{o}}$lder inequality and using estimates for $\mathcal{\widetilde U}_\ep$ in (\ref{estim_u_I_dil})--(\ref{estim_Du_I_dil}), we  get
$$
\begin{array}{rl}
\displaystyle
\left|-\nu\int_{\widetilde I_1}|\mathbb{D}_{\eta_\ep}[\mathcal{\widetilde U}_\ep]|^{r-2}\mathbb{D}_{\eta_\ep}[\mathcal{\widetilde U}_\ep]\,dz\right|\leq &\displaystyle C\|\mathbb{D}_{\eta_\ep}[\mathcal{\widetilde U}_\ep]\|_{L^r(\widetilde I_1)^{2\times 2}}^{r-1}\|D_{\eta_\ep} \widetilde \varphi\|_{L^r(\widetilde I_1)^{2\times 2}}\\
 \leq& C\left( \eta_\ep+h_\ep^{r-1\over r}\ep\eta_\ep^{-{r-1\over r}}\right)\|D_{\eta_\ep}\widetilde \varphi\|_{L^r(\widetilde I_1)^{2\times 2}}\\
 \leq&\displaystyle C\left( 1+h_\ep^{r-1\over r}\ep\eta_\ep^{-{2r-1\over r}}\right)\|\widetilde \varphi\|_{W^{1,r}_0(\widetilde I_1)^2},
\end{array}
$$
where, in the last inequality, we have used  
\begin{equation}\label{TESTPHI}\|D_{\eta_\ep}\widetilde\varphi\|_{L^r(\widetilde I_1)^{2\times 2}}\leq C\eta_\ep^{-1}\|\widetilde \varphi\|_{W^{1,r}_0(\widetilde I_1)^2}.
\end{equation}
\item[--] Applying the change of variables (\ref{dilatacion_t}) to inequality (\ref{PoincareIep1}), we get
$$
\|\widetilde \varphi\|_{L^r(\widetilde I_1)^2}\leq C\eta_\ep^{1\over 2}(\ep+\eta_\ep)^{1\over 2}\|D_{\eta_\ep}\widetilde \varphi\|_{L^r(\widetilde I_1)^{2\times 2}},
$$
and from the H${\rm \ddot{o}}$lder inequality and using this inequality, assumption of $f$ given in (\ref{fep})   and (\ref{TESTPHI}), we have
$$\begin{array}{rl}
\displaystyle
 \left|\int_{\widetilde I_1}f\cdot \widetilde \varphi \,dz\right|\leq&\displaystyle  C\| \widetilde\varphi \|_{L^r(\widetilde I_1)^2}\leq  C\eta_\ep^{1\over 2}(\eta_\ep+\ep)^{1\over 2}\|D_{\eta_\ep}\widetilde \varphi\|_{L^r(\widetilde I_1)^{2\times 2}}\\
 \leq &\displaystyle C\eta_\ep^{-{1\over 2}}(\eta_\ep+\ep)^{1\over 2}\|\widetilde \varphi\|_{W^{1,r}_0(\widetilde I_1)^2}\leq C\left(1+  \eta_\ep^{-{1\over 2}}\ep^{1\over 2} \right)\|\widetilde \varphi\|_{W^{1,r}_0(\widetilde I_1)^2}\,.
\end{array}$$
\end{itemize}
Previous estimates together with (\ref{equality_duality23}) gives $$\left|\left\langle \nabla_{\eta_\ep}\widetilde P^2_\ep,\widetilde \varphi\right\rangle_{W^{-1,r'}(\widetilde I_1)^2, W^{1,r}_0(\widetilde I_1)^2}\right|\leq
 C\left( 1+h_\ep^{r-1\over r}\ep\eta_\ep^{-{2r-1\over r}}+\eta_\ep^{-{1\over 2}}\ep^{1\over 2}
 \right)\|\widetilde \varphi\|_{W^{1,r}_0(\widetilde I_1)^2}.$$
Since $\eta_\ep\gg\ep$, then the term $\eta_\ep^{-{1\over 2}}\ep^{1\over 2}\ll 1$ and it can be dropped. This implies the second estimate given in (\ref{esti_P_thin_film}), which concludes the proof.
\end{proof}

\begin{remark}
In view of estimates of the velocity and the pressure given in the previous section,
there is a critical case when 
\begin{equation}\label{CritCase}
h_\ep\approx  \eta_\ep^{2r-1\over r-1} \ep^{-{r\over r-1}}\quad\hbox{with}\quad \lim_{\ep\to 0}{h_\ep\over  \eta_\ep^{2r-1\over r-1} \ep^{-{r\over r-1}}}= \lambda,\quad 0<\lambda<+\infty,
\end{equation} 
where the pressure has the same order of magnitude in the porous medium and in the free film. From now on, we focus our study in this case, which is the most interesting one.
\end{remark}

\section{Critical case: problem in the thin porous medium}\label{sec:thinporousmedium}

In this section, we study the asymptotic behavior of the fluid in the thin porous part assuming the critical regime (\ref{CritCase}).  Under this assumption, the estimates given in Corollary \ref{estimates_vel_tilde} and Lemma \ref{Estimates_extended_lemma} read as follow
\begin{equation}\label{estimates_critical_case}
\begin{array}{c}
\displaystyle
\|\widetilde v_\ep\|_{L^r(\widetilde \Omega_\ep)^{2}}\leq C\ep^{r\over r-1},\quad \|\mathbb{D}_{h_\ep}[\widetilde v_\ep]\|_{L^r(\widetilde \Omega_\ep)^{2\times 2}}\leq C\ep^{1\over r-1},\quad \|D_{h_\ep}\widetilde v_\ep\|_{L^r(\widetilde \Omega_\ep)^{2\times 2}}\leq C\ep^{1\over r-1},\\
\\ \|\widetilde P_\ep^1\|_{L^{r'}(\Omega)}\leq C,\quad \|\nabla_{h_\ep}  \widetilde P_\ep^1\|_{W^{-1,r'}(\Omega)^2}\leq C.
\end{array}
\end{equation}
To describe the behavior of the solution in the microstructure associated to $\widetilde\Omega_\ep$, we introduce an adaptation of the unfolding method (for classical versions see  \cite{Ciora2,Cioran-book}), which is related with the change of variables applied in \cite{BayadaThinThin} (see also \cite{SG_MANA}) to study the porous part in the Newtonian case of modeling of a thin film passing a thin porous media.

\subsection{Adaptation of the unfolding method}\label{sec:unfolding}
This version of the unfolding method  consists of dividing the domain $\widetilde\Omega_\ep$ into squares of horizontal length $\ep$ and vertical length $\ep/h_\ep$. In order to apply the version of the unfolding method, we need the following notation: for $k\in\mathbb{Z}^2$, we define $\kappa:\mathbb{R}^2\to \mathbb{Z}^2$ by
\begin{equation}\label{kappa}
\kappa(x)=k\Longleftrightarrow x\in Y_{k,1}.
\end{equation}
Remark that $\kappa$ is well defined up to a set of zero measure in $\mathbb{R}^2$, which is given by $\cup_{k\in\mathbb{Z}^2}\partial Y_{k,1}$.  Moreover, for every $\ep,\,h_\ep>0$, we have 
$$\kappa\left({x\over \ep}\right)=k\Longleftrightarrow x\in Y_{k,\ep}\quad\hbox{which is equivalent to }\quad \kappa\left({z_1\over \ep},{h_\ep z_2\over \ep}\right)=k\Longleftrightarrow z\in \widetilde Y_{k,\ep}.$$

\begin{definition} Let $\widetilde \varphi$ be in $L^{s}(\widetilde\Omega_\ep)^2$, $1\leq s<+\infty$,  and $\widetilde \psi$ be in $L^{s'}(\Omega)$, $1/s+1/s'=1$. We define the functions $\widehat \varphi_\ep\in L^s(\mathbb{R}^2\times Y_f)^2$ and $\widehat \psi_\ep\in L^{s'}(\mathbb{R}^2\times Y)$ by 
\begin{eqnarray}
\displaystyle\widehat \varphi_\ep(z,y)= \widetilde \varphi\left( {\varepsilon}\kappa\left(\frac{z_1}{{\varepsilon}}, { h_\ep z_2\over \ep } \right)\cdot e_1+ {\varepsilon}y_1, {{\varepsilon}\over h_\ep}\kappa\left(\frac{z_1}{ {\varepsilon}}, { h_\ep z_2\over  \ep   } \right)\cdot e_2+{ {\varepsilon}\over h_\ep}y_2 \right),& \hbox{a.e. }(z, y)\in  \mathbb{R}^2 \times   Y_f,&\label{def:unfolding1}\\\nonumber
\\
\displaystyle
\widehat \psi_\ep(z,y)=  \widetilde \psi\left( {\varepsilon}\kappa\left(\frac{z_1}{{\varepsilon}}, { h_\ep z_2\over \ep } \right)\cdot e_1+ {\varepsilon}y_1, {{\varepsilon}\over h_\ep}\kappa\left(\frac{z_1}{ {\varepsilon}}, { h_\ep z_2\over  \ep   } \right)\cdot e_2+{ {\varepsilon}\over h_\ep}y_2 \right),&  \hbox{a.e. }(z, y)\in \mathbb{R}^2\times   Y,&\label{def:unfolding2}
\end{eqnarray}
assuming $\widetilde \varphi$ (resp. $\widetilde \psi$) is extended by zero outside $\widetilde\Omega_\ep$ (resp. $\Omega$), where
 the function $\kappa$ is defined by (\ref{kappa}).
\end{definition}
\begin{remark}\label{remarkCV}
The restrictions of $\widehat \varphi_{\varepsilon}$  to $\widetilde Y_{k,{\varepsilon}}\times Y_f$  (resp.  $\widehat \psi_\ep$ to $\widetilde Y_{k,{\varepsilon}}\times Y$) does not depend on $z$, while as a function of $y$ it is obtained from $\widetilde \varphi$ (resp.  from $\widetilde \psi)$ by using the change of variables 
\begin{equation}\label{CV}
y_1=\frac{z_1-{\varepsilon}k_1}{{\varepsilon}},\quad y_2=\frac{h_\ep z_2-{\varepsilon}k_2}{{\varepsilon}},
\end{equation}
which transforms $\widetilde Y_{f_k,{\varepsilon}}$ into $Y_f$ (resp. $\widetilde Y_{k,{\varepsilon}}$ into $Y$).
\end{remark}
Next, we give some properties of the unfolded functions.
\begin{proposition} \label{properties_um}Let $\widetilde \varphi$ be in $W^{1,s}(\widetilde\Omega_\ep)^2$, $1\leq s<+\infty$,  and $\widetilde \psi$ be in $L^{s'}(\Omega)$, $1/s+1/s'=1$. Then, we have
\begin{equation*}\label{relation_norms}
\begin{array}{lcl}
\displaystyle
\|\widehat \varphi_\ep\|_{L^s(\mathbb{R}^2\times Y_f)^2}=\|\widetilde \varphi\|_{L^s(\widetilde \Omega_\ep)^2},& \|\partial_{y_1}\widehat \varphi_\ep\|_{L^s(\mathbb{R}^2\times Y_f)^{2}}=\ep\|\partial_{z_1}\widetilde \varphi\|_{L^s(\widetilde \Omega_\ep)^{2}},&\displaystyle \|\partial_{y_2}\widehat \varphi_\ep\|_{L^s(\mathbb{R}^2\times Y_f)^2}={\ep\over h_\ep}\|\partial_{z_2}\widetilde \varphi\|_{L^s(\widetilde \Omega_\ep)^2},\\
\\
& 
\|\widehat \psi_\ep\|_{L^{r'}(\mathbb{R}^2\times Y)}=\|\widetilde \psi\|_{L^{r'}(\Omega)}.&

\end{array}\end{equation*}
\end{proposition}
\begin{proof} We will only make the proof for $\widehat \varphi_\ep$. The procedure for $\widehat \psi_\varepsilon$ is similar, so we omit it. Taking into account the definition (\ref{def:unfolding1}) of $\widehat{\varphi}_{\varepsilon}$, we obtain
\begin{eqnarray*}
\int_{\mathbb{R}^2\times Y_f}\left\vert \partial_{y_1} \widehat{\varphi}_{\varepsilon}(z, y) \right\vert^s dzdy&=&\displaystyle\sum_{k\in  \mathbb{Z}^2}\int_{\widetilde Y_{k,{\varepsilon}}}\int_{Y_f}\left\vert \partial_{y_1} \widehat{\varphi}_{\varepsilon}(z, y) \right\vert^s dy dz\\
&=&\displaystyle\sum_{k \in \mathbb{Z}^2}\int_{\widetilde Y_{k,{\varepsilon}}}\int_{Y_f}\left\vert \partial_{y_1} \widetilde{\varphi} ({\varepsilon}k^{\prime}+{\varepsilon}y_1,{\varepsilon}h_\ep^{-1}k_2+ {\varepsilon}h_\ep^{-1}y_2
) \right\vert^sdydz.
\end{eqnarray*}
We observe that $\widetilde{\varphi}$ does not depend on $z$, then we can deduce
\begin{eqnarray*}
\int_{\mathbb{R}^2\times Y_f}\left\vert \partial_{y_1} \widehat{\varphi}_{\varepsilon}(z, y) \right\vert^s dzdy
= {{\varepsilon}^2\over h_\ep}\displaystyle\sum_{k \in \mathbb{Z}^2}\int_{Y_f}\left\vert \partial_{y_1} \widetilde{\varphi}({\varepsilon}k_1+ {\varepsilon}y_1,{\varepsilon}h_\ep^{-1}k_2+{\varepsilon}h_\ep^{-1}y_2
) \right\vert^sdy.
\end{eqnarray*}
By the change of variables (\ref{CV}), we obtain
$$
\int_{\mathbb{R}^2\times Y_f}\left\vert \partial_{y_1} \widehat{\varphi}_{\varepsilon}(z, y) \right\vert^s dzdy
={\varepsilon}^s
\displaystyle\sum_{k \in  \mathbb{Z}^2}\int_{\widetilde Y_{f_k,{\varepsilon}}} \left\vert \partial_{y_1} \widetilde{\varphi}(z) \right\vert^sdz= {\varepsilon}^s\int_{\widetilde \Omega_\ep}\left\vert \partial_{x_1} \widetilde \varphi(z) \right\vert^sdz.
$$
Thus, we get the property for $\partial_{y_1}\widehat \varphi_\ep$.

Similarly,  we have
\begin{eqnarray*}
\int_{\mathbb{R}^2\times Y_f}\left\vert \partial_{y_2} \widehat{\varphi}_{\varepsilon}(z, y) \right\vert^sdzdy= {{\varepsilon}^2\over h_\ep}\displaystyle\sum_{k \in \mathbb{Z}^2}\int_{Y_f}\left\vert \partial_{y_2} \widetilde{\varphi} ({\varepsilon}k_1+{\varepsilon}y_1, {\varepsilon}h_\ep^{-1}k_2+{\varepsilon}h_\ep^{-1}y_2
)\right\vert^sdy.
\end{eqnarray*}
By the change of variables (\ref{CV})  we obtain
$$
\int_{\mathbb{R}^2\times Y_f}\left\vert \partial_{y_2} \widehat{\varphi}_{\varepsilon}(z,y) \right\vert^sdzdy 
=
{{\varepsilon}^s\over h_\ep^s}\displaystyle\sum_{k \in \mathbb{Z}^2}\int_{\widetilde Y_{f_k,\ep}}\left\vert \partial_{z_2} \widetilde{\varphi} (z)\right\vert^sdz={{\varepsilon}^s\over h_\ep^s} \int_{\widetilde\Omega_\ep}\left\vert \partial_{z_2} \widetilde{\varphi}(z)\right\vert^sdz,
$$
so the the property for $\partial_{y_2}\widehat \varphi_\ep$ is proved. Finally, reasoning analogously we deduce
\begin{eqnarray*}
\int_{\mathbb{R}^2\times Y_f}\left\vert \widehat{\varphi}_{\varepsilon}(z,y)\right\vert^sdzdy = \int_{\widetilde\Omega_\ep}\left\vert \widetilde{\varphi}(z)\right\vert^sdz,
\end{eqnarray*}
and the property for $\widehat \varphi_\ep$  holds. 
\end{proof}

\begin{lemma}\label{estimates_hat} We assume that the parameters $\ep, \eta_\ep$ and $h_\ep$ satisfy (\ref{parameters}) and (\ref{CritCase}). We define the unfolded velocity $\widehat v_\ep$  from  the dilated velocity $\widetilde v_\ep$ by means of (\ref{def:unfolding1}) and the unfolded pressure $\widehat P_\ep^1$ from the dilated and extended pressure $\widetilde P_\ep^1$ by means of (\ref{def:unfolding2}). Then, there exists a constant $C>0$ independent of $\ep$, such that $\widehat v_\ep$   and   $\widehat P_\ep^1$  satisfy  
\begin{equation}\label{estim_u_hat}
 \|\widehat v_\ep\|_{L^r(\mathbb{R}^2\times Y_f)^2}\leq C\ep^{r\over r-1},\quad 
 \|D_{y}\widehat v_\ep\|_{L^r(\mathbb{R}^2\times Y_f)^{2\times 2}}\leq C\ep^{r\over r-1},
 \end{equation}
\begin{equation}\label{estim_P_hat}
 \|\widehat P_\ep^1\|_{L^{r'}(\mathbb{R}^2\times Y)}\leq C.
\end{equation}
 \end{lemma}
\begin{proof} 
Estimates (\ref{estim_u_hat}) and (\ref{estim_P_hat}) easily follow from Proposition \ref{properties_um}, with $s=r$ and $s'=r'$, and  estimates of velocity in $\widetilde\Omega_\ep$ and estimate of the pressure in $\Omega$ given in (\ref{estimates_critical_case}).
\end{proof}

\subsection{Convergences of velocity and pressure}\label{sec:compactness}
From now on, we denote by $\widetilde V_\ep$ the extension by zero of $\widetilde v_\ep$ to the whole domain $\Omega$ (the velocity is zero in the obstacles). Then, estimates given in (\ref{estimates_critical_case}) remain valid for the extension $\widetilde V_\ep$, which is divergence free too.  Here, we obtain some compactness results concerning the behavior of the sequence $(\widetilde V_\ep, \widetilde P_\ep^1)$ and $(\widehat v_\ep,\widehat P_\ep^1)$.
\begin{lemma}\label{lemma_compactness}
For a subsequence of $\ep$ still denoted by $\ep$,  there exist:
\begin{itemize}

\item[--]    $v\in W^{1,r}(0,1;L^r(\omega)^2)$, with $v_2\equiv 0$ and $v_1=0$ on $\{z_2=1\}$, such that
\begin{eqnarray}
&\ep^{-{r\over r-1}}\widetilde V_\ep \rightharpoonup v\quad\hbox{in }W^{1,r}(0,1;L^r(\omega)^2),&\label{conv_vel_tilde}
\end{eqnarray}
\begin{equation}\label{divxproperty}
{\partial}_{z_1} \left(\int_{0}^1v_1(z)\,dz_2\right)=0\quad \hbox{in }\omega,\quad   \left(\int_{0}^1v_1(z)\,dz_2\right)\, n=0\quad \hbox{on }\partial\omega.
\end{equation}
\item[--]  $\widehat v\in L^r(\mathbb{R}^2; W^{1,r}_{{\rm per}}(Y_f)^2)$, with  $\widehat v_2$ independent of $z_2$ and $\widehat v=0$ in $\Omega\times  T$ and in $(\mathbb{R}^2\setminus \Omega)\times Y_f$, such that
\begin{eqnarray}
&\ep^{-{r\over r-1}}\widehat v_\ep\rightharpoonup \widehat v\hbox{ in }L^r(\mathbb{R}^2; W^{1,r}(Y_f)^2),\quad \ep^{-{r\over r-1}}D_y\widehat v_\ep\rightharpoonup D_y\widehat v\hbox{ in }L^r(\mathbb{R}^2\times Y_f)^{2\times 2}.\label{conv_vel_gorro}&
\end{eqnarray}
\begin{eqnarray}
& \displaystyle {\rm div}_{y}\,\widehat v(z,y)=0\quad \hbox{in }\mathbb{R}^2\times Y_f.&\label{divyproperty}
\end{eqnarray}
\end{itemize}
Moreover, the following relation between $v$ and $\widehat v$ holds 
 \begin{equation}\label{relation_u_ugorro}
 v(z)= \int_{Y_f}\widehat v(z,y)\,dy\quad  \hbox{a.e. in  }\Omega,\hbox{ i.e.  }\quad\int_{Y_f}\widehat v(z,y)\,dy=v_1(z)\quad \hbox{and}\quad \int_{Y_f}\widehat v_2(z,y)\,dy=0.
 \end{equation}
\end{lemma}
\begin{proof}We divide the proof in two parts:
\begin{itemize} 
\item[--] We start with the extended velocity $\widetilde V_\ep$. From the first and second estimate in (\ref{estimates_critical_case}), 
 we get the existence of $v\in W^{1,r}(0,1;L^r(\omega)^2)$ such that, up to a subsequence, it holds
\begin{equation}\label{convergencia_debil_utilde}
\ep^{-{r\over r-1}}\widetilde V_\ep \rightharpoonup v\hbox{ in }W^{1,r}(0,1;L^r(\omega)^2).
\end{equation}
This implies
\begin{equation}\label{convergencia_debil_utilde_dual}
\ep^{-{r\over r-1}}{\partial}_{z_1}\widetilde V_{\ep,1} \rightharpoonup {\partial}_{z_1}v_1\hbox{ in }W^{1,r}(0,1;W^{-1,r'}(\omega)^2).
\end{equation} 
Since ${\rm div}_{h_\ep}(\widetilde V_\ep)=0$ in $\Omega$, multiplying by $h_\ep \ep^{-{r\over r-1}}$ we obtain
\begin{equation}\label{div_limit1}
h_\ep \ep^{-{r\over r-1}}{\partial}_{z_1}\widetilde V_{\ep,1}+ \ep^{-{r\over r-1}}\partial_{z_2} \widetilde V_{\ep,2}=0\quad\hbox{in }\Omega,
\end{equation}
which, combined with (\ref{convergencia_debil_utilde_dual}), implies that $\ep^{-{r\over r-1}}\partial_{z_2} \widetilde V_{\ep,2}$ is bounded in $W^{1,r}(0,1;W^{-1,r'}(\omega)^2)$ and tends to zero. Also from (\ref{convergencia_debil_utilde}), we have that $\ep^{-{r\over r-1}}\partial_{z_2} \widetilde V_{\ep,2}$ tends to $\partial_{z_2} v_{2}$ in $L^r(\Omega)$. From the uniqueness of the limit, we have that $\partial_{z_2} v_{2}=0$, which implies that $v_2$ is independent of $z_2$. Moreover, the continuity of the trace applications from the space of functions $\widetilde \varphi$ such that $\|\widetilde \varphi\|_{L^r}$ and $\|\partial_{z_2}\widetilde \varphi\|_{L^r}$ to $L^r(\omega\times \{1\})$ implies $v=0$ on $z_2=\{1\}$. From this boundary condition and since $v_2$ does not depend on $z_2$, we deduce $v_2\equiv 0$. This completes the proof of (\ref{conv_vel_tilde}).

 Next, by considering $\widetilde \varphi\in \mathcal{D}(\omega)$ as test function in the divergence condition ${\rm div}_{h_\ep}\widetilde V_\ep=0$ in $\Omega$, we get
$$\int_{\Omega}\left({\partial}_{z_1}\widetilde V_{\ep,1}\, \widetilde \varphi+h_\ep^{-1}\partial_{z_2}\widetilde V_{\ep,2}\widetilde  \varphi\right)\,dz=0,$$
which, after integration by parts and multiplication by $\ep^{-{r\over r-1}}$, gives
$$\int_{\Omega} \ep^{-{r\over r-1}}\widetilde V_{\ep,1}\, \partial_{z_1}\widetilde \varphi\,dz=0.$$
Passing to the limit by using convergence (\ref{conv_vel_tilde}), we deduce
$$\int_{\Omega}v_1\, \partial_{z_1}\widetilde \varphi\,dz=0,$$
and, since $\widetilde \varphi$ does not depend on $z_2$, we obtain the following divergence condition (\ref{divxproperty}).

\item[--] Now we focus on the velocity $\widehat v_\ep$. From estimates of  $\widehat v_\ep$ given in (\ref{estim_u_hat}) we have the existence of $\widehat v\in L^r(\mathbb{R}^2; W^{1,r}_{\rm per}(Y_f)^2)$ satisfying, up to a subsequence,  convergences (\ref{conv_vel_gorro}). Taking into account that $\widehat v_\ep$ vanishes on $\mathbb{R}^2\times T$, we deduce that $\widehat v$ also vanishes on  $\mathbb{R}^2\times T$.  Moreover, by construction $\widehat v_\ep$ is zero outside $\widetilde\Omega_\ep$ and so, $\widehat v$ vanishes on $(\mathbb{R}^2\setminus \Omega)\times Y_f$.

Since ${\rm div}_{h_\ep}(\widetilde v_\ep)=0$ in $\widetilde \Omega_\ep$,  by applying the change of variables (\ref{CV}) we get 
$$\ep^{-1}{\rm div}_y(\widehat v_\ep)=0\quad \hbox{in }\mathbb{R}^2\times Y_f.$$
Multiplying by $\ep^{-{1\over r-1}}$ and passing to the limit by using convergence (\ref{conv_vel_gorro}), we deduce ${\rm div}_y(\widehat v)=0$ in $\mathbb{R}^2\times Y_f$, i.e. property (\ref{divyproperty}).
\\

It remains to prove that $\widehat v$ is periodic in $y$. This follows by passing to the limit in the equality
$$\ep^{-{r\over r-1}}\widehat v_\ep\left(z+\ep\,{\rm e}_1,-{1\over 2}, y_2\right)=\ep^{-{r\over r-1}}\widehat v_\ep\left(z,{1\over 2}, y_2\right),$$
which is a consequence of definition (\ref{def:unfolding1}). Passing to the limit, this shows 
$$\widehat v\left(z,-{1\over 2},y_2\right)=\widehat v\left(z,{1\over 2},y_2\right),$$
and then is proved  the periodicity of $\widehat v$ with respect to $y_1$. To prove the periodicity with respect to $y_2$, we consider 
$$\ep^{-{r\over r-1}}\widehat v_\ep\left( z+{\ep\over h_\ep}{\rm e}_2,y_1,-{1\over 2}\right)=\ep^{-{r\over r-1}}\widehat v_\ep\left(z,y_1,{1\over 2}\right),$$
and passing to the limit we have 
$$\widehat v\left(z,y_1,-{1\over 2}\right)=\widehat v\left(z,y_1,{1\over 2}\right),$$
which shows the periodicity with respect to $y_2$.

Finally, relation (\ref{relation_u_ugorro}) follows from Proposition \ref{properties_um} with $s=1$, which gives
$$\int_\Omega v(z) \,dz=\int_{\Omega\times Y_f}\widehat v(z,y)\,dzdy=\int_{\Omega}\left(\int_{Y_f}\widehat v(z,y)\,dy\right)\,dz.$$
From relation (\ref{relation_u_ugorro}) and since $v_2\equiv 0$, it holds that $\int_{Y_f}\widehat v_2\,dy=0$.

\end{itemize}
\end{proof}

\begin{lemma}\label{lemma_conv_pressure}
For a subsequence of $\ep$ still denoted by $\ep$, there exists  $p^1\in L^{r'}_0(\omega)$ independent of $z_2$, such that 
 \begin{equation}\label{conv_pressure_sub}
\widetilde P_\ep^1\to p^1\quad\hbox{ in }L^{r'}(\Omega),
\end{equation} 
 \begin{equation}\label{conv_pressure_gorro}
\widehat P_\ep^1 \to p^1 \quad\hbox{ in }L^{r'}(\mathbb{R}^2\times Y).
\end{equation}

\end{lemma}
\begin{proof} Taking into account the first estimate of the pressure in (\ref{estimates_critical_case}), we deduce that there exist $p^1\in L^{r'}(\Omega)$ such that, up to a subsequence, 
\begin{equation}\label{conv_pressure_sub_weak}
\widetilde P_\ep^1\rightharpoonup p^1\quad\hbox{ in }L^{r'}(\Omega).
\end{equation} 
 From convergence (\ref{conv_pressure_sub_weak}) we deduce that $\partial_{z_2}\widetilde P_\ep^1$ also converges to $\partial_{z_2}p^1$ in $W^{-1,r'}(\Omega)$. Also, from the second estimate of the pressure in (\ref{estimates_critical_case}), we can deduce that  $\partial_{z_2}\widetilde P_\ep^1$ converges to zero in  $W^{-1,r'}(\Omega)$.  By the uniqueness of the limit, then  we obtain $\partial_{z_2} p^1=0$ and so $p^1$ is independent of $z_2$. Since $\widetilde P_\ep^1$ has null mean value in $\Omega$, then $p^1$ has null mean value in $\omega$.

Next, following  \cite{Bourgeat1} adapted to the case of a thin layer, we prove that the convergence of the pressure is in fact strong. 
 Let $\widetilde w_\ep, \widetilde w$ be in  $W^{1,r}_0(\Omega)^2$  such that
\begin{equation}\label{strong_p_1}
\widetilde w_\ep\rightharpoonup \widetilde w\quad\hbox{in }W^{1,r}_0(\Omega)^2.
\end{equation}
Then, as $p^1$ only depends on $z_1$, we have 
$$\begin{array}{l}
\displaystyle
\left|\langle\nabla_{z}\widetilde P_\ep^1,\widetilde w_\ep\rangle_{W^{-1,r'}(\Omega)^2,W^{1,r}_0(\Omega)^2}-\langle\nabla_{z}p^1,\widetilde w\rangle_{W^{-1,r'}(\Omega)^2,W^{1,r}_0(\Omega)^2}\right| 
\\
\\
 \displaystyle \leq 
\left|\langle\nabla_{z}\widetilde P_\ep^1,\widetilde w_\ep-\widetilde w\rangle_{W^{-1,r'}(\Omega)^2,W^{1,r}_0(\Omega)^2}\right| +\left|\langle\nabla_{z}(\widetilde P_\ep^1-p^1),\widetilde w\rangle_{W^{-1,r'}(\Omega)^2,W^{1,r}_0(\Omega)^2}\right|.
\end{array}$$
On the one hand, using convergence (\ref{conv_pressure_sub_weak}), we have 
$$\left|\langle\nabla_{z} (\widetilde P_\ep^1- p^1),\widetilde w\rangle_{W^{-1,r'}(\Omega)^2,W^{1,r}_0(\Omega)^2}\right|=\left|\int_\Omega\left(\widetilde P_\ep^1-p^1\right)\,{\rm div}_z\widetilde w\,dz\right|\to 0,\quad \hbox{as }\ep\to 0\,.$$
On the other hand, from (\ref{equality_duality2}) and proceeding similarly to the proof of Lemma \ref{Estimates_extended_lemma}, we have
$$\begin{array}{rl}
\left|\langle\nabla_{z}\widetilde P_\ep^1,\widetilde w_\ep-\widetilde w\rangle_{W^{-1,r'}(\Omega)^2,W^{1,r}_0(\Omega)^2}\right|\leq & \left|\langle\nabla_{h_\ep}\widetilde P_\ep^1,\widetilde w_\ep-\widetilde w\rangle_{W^{-1,r'}(\Omega)^2,W^{1,r}_0(\Omega)^2}\right|\\
\\\displaystyle
 \leq &C\left(\|\widetilde w_\ep-\widetilde w\|_{L^r(\Omega)^2}+\ep\|D_{h_\varepsilon} (\widetilde w_\ep-\widetilde w)\|_{L^r(\Omega)^{2\times 2}}\right)
 \\
\\\displaystyle
 \leq &C\left(\|\widetilde w_\ep-\widetilde w\|_{L^r(\Omega)^2}+\ep h_\ep^{-1}\|D_{z} (\widetilde w_\ep-\widetilde w)\|_{L^r(\Omega)^{2\times 2}}\right).
\end{array}$$
The right-hand side of the previous inequality tends to zero as $\ep\to 0$, by virtue of relation (\ref{parameters}), (\ref{strong_p_1}) and the Rellich theorem. This implies that $\nabla_{z}\widetilde P_\ep^1\to \nabla_z p^1=(\partial_{z_1} p^1,0)^t$ strongly in $W^{-1,r'}(\Omega)^3$, which together the classical Ne${\breve{\rm c}}$as inequality implies the strong convergence of the pressure $\widetilde P_\ep^1$ given in (\ref{conv_pressure_sub}).  Finally,   the strong convergence of $\widehat P_\varepsilon^1$ given in (\ref{conv_pressure_gorro}) follows from \cite[Proposition 1.9-(ii)]{Cioran-book} and the strong convergence of $\widetilde P_\varepsilon^1$ given in (\ref{conv_pressure_sub}).

\end{proof}
\subsection{Average velocity in the porous medium}
We deduce an expression for the average limit velocity in the thin porous medium.
\begin{theorem}\label{mainthm_porous}Consider the pair of limit functions $(\widehat v, p^1)$  given in Lemmas \ref{lemma_compactness} and \ref{lemma_conv_pressure}. Defining the  average velocity  by 
$$V_{av}(z_1)=\int_0^1\int_{Y_f}\widehat v(z,y)\,dydz_2,$$  we have 
\begin{equation}\label{VelProm}V_{av,1}(z_1)={\mu\over \nu^{r'-1}}\Big|f_1(z_1)-{d\over d{z_1}}p^1(z_1)\Big|^{r'-2}\Big(f_1(z_1)-{d\over d{z_1}}p^1(z_1)\Big),\quad V_{av, 2}\equiv 0,\quad \hbox{in }\omega,
\end{equation}
and from (\ref{divxproperty}) and $\int_0^1v_1(z)dz_2=V_{av,1}(z_1)$, we have
\begin{equation}\label{property_pressure}
\partial_{z_1}(V_{av,1}(z_1))=0,\quad\hbox{in }\omega,\qquad V_{av,1}\cdot n=0\quad\hbox{on }\partial\omega,
\end{equation}
where $\mu\in\mathbb{R}$ is the permeability defined by 
\begin{equation}\label{nueff}
\mu=\int_{Y_f}\widehat w(y)\cdot {\rm e}_1\,dy,
\end{equation}
where $(\widehat w, \widehat q)\in W^{1,r}_{\rm per}(Y)^2\times L^{r'}_{0, {\rm per}}(Y)$, $1<r<+\infty$, is the unique solution of the auxiliary problem 
\begin{equation}\label{Local_problems}
\left\{\begin{array}{rl}
-{\rm div}_y\left(|\mathbb{D}_y[\widehat w]|^{r-2}\mathbb{D}_y[\widehat w]\right)+ \nabla_y \widehat q={\rm e}_1 & \hbox{ in }Y_f,\\
\\
{\rm div}_y \widehat w=0& \hbox{ in }Y_f,\\
\\
\widehat w=0& \hbox{ on }T.
\end{array}\right.
\end{equation}
\end{theorem}
\begin{remark}
As is pointed in \cite[Remark 8]{Bourgeat1}, we observe that we have derived a Darcy law (\ref{VelProm}) identical the usual filtration law used in standard engineering treatment (see for instance Wu {\it et al} \cite[p. 140]{Wu}).   We point out that the version of the unfolding method and the restriction operator introduced in this paper are powerfull tools that could be used to study the asymptotic behavior of different type of (two dimensional or three dimensional) fluids in a thin porous medium defined by $\Omega_\ep$. 
\end{remark}
\begin{proof}[Proof of Theorem \ref{mainthm_porous}]  We divide the proof in three steps.

{\it Step 1. Variational formulation for $(\widehat v_\ep,\widehat P_\ep^1)$.} Let us first write the variational formulation satisfied by the functions $(\widehat u_\ep,\widehat P_\ep^1)$ in order to pass to the limit. According to Lemma \ref{lemma_compactness}, we consider $\widetilde \varphi_\ep(z)=(\widehat \varphi_1(z_1,  z_2,z_1/\ep,$ $h_\ep z_2/\ep),$ $ \widehat \varphi_2(z_1,z_1/\ep,h_\ep z_2/\ep))$,  as test function in (\ref{form_var_1_tildev}) where $\widehat \varphi(z,y)=(\widehat \varphi_1(z,y), \widehat \varphi_2(z_1,y))  \in \mathcal{D}(\Omega;C^\infty_{\rm per}(Y)^2)$ with $\widehat \varphi(z,y)=0$ in $\Omega\times T$ and $(\mathbb{R}^2\setminus \Omega)\times Y$ (thus, $\widetilde \varphi_\ep(z)\in W^{1,r}_0(\widetilde\Omega_\ep)^2$). Then, we have
\begin{equation}\label{form_var_general_a}
\begin{array}{l}\displaystyle
\nu\int_{\widetilde\Omega_\ep}S_r(\mathbb{D}_{h_\ep}[\widetilde v_\ep]):\mathbb{D}_{h_\ep}\widetilde \varphi_\ep\,dz+\langle \nabla_{h_\ep} \widetilde p_\ep,\widetilde\varphi_\ep\rangle_{W^{-1,r'}(\widetilde\Omega_\ep)^2, W^{1,r}_0(\widetilde\Omega_\ep)^2} 
=\int_{\widetilde\Omega_\ep}f_1\, (\widetilde \varphi_\ep)_1\,dz\,,
\end{array}
\end{equation}
where, for simplicity, we have denoted by $S_r: \mathbb{R}^{2\times 2}_{\rm sym}\to \mathbb{R}^{2\times 2}_{\rm sym}$   the $r$-Laplace operator, i.e. $S_r$ is defined by
\begin{equation}\label{opS}
S_r(\xi)=|\xi|^{r-2}\xi,\quad \forall\,\xi\in \mathbb{R}^{2\times 2}_{\rm sym},\quad 1<r<+\infty.
\end{equation}
Taking into account the extension of the pressure, we get
$$\langle\nabla_{h_\ep} \widetilde p_\ep,  \widetilde \varphi_\ep\rangle_{W^{-1,r'}(\widetilde\Omega_\ep)^2, W^{1,r}_0(\widetilde\Omega_\ep)^2}  =\langle\nabla_{h_\ep} \widetilde P_\ep^1,  \widetilde \varphi_\ep\rangle_{W^{-1,r'}(\Omega)^2, W^{1,r}_0(\Omega)^2} =- \int_{\Omega}  \widetilde P_\ep^1\,{\rm div}_{h_\ep}(\widetilde \varphi_\ep)\,dz,$$
and then, the variational formulation (\ref{form_var_general_a}) reads
\begin{equation}\label{form_var_general_1}
\begin{array}{l}\displaystyle
\nu\int_{\widetilde\Omega_\ep}S_r(\mathbb{D}_{h_\ep}[\widetilde v_\ep]):\mathbb{D}_{h_\ep}\widetilde \varphi_\ep\,dz-\int_{\Omega}  \widetilde P_\ep^1\,{\rm div}_{h_\ep}(\widetilde \varphi_\ep)\,dz
=\int_{\widetilde\Omega_\ep}f_1\, (\widetilde \varphi_\ep)_1\,dz\,.
\end{array}
\end{equation}
Taking into account the definition of $\widetilde \varphi_{\varepsilon}$, we have
$$\begin{array}{l}
\displaystyle 
\partial_{z_1}\widetilde\varphi_{\varepsilon,1}(z)=\partial_{z_1}\widehat \varphi_1+ \varepsilon^{-1}\partial_{y_1}\widehat \varphi_1,\quad 
h_\varepsilon^{-1}\partial_{z_1}\widetilde\varphi_{\varepsilon,1}(z)=h_\varepsilon^{-1}\partial_{z_1}\widehat \varphi_1+ \varepsilon^{-1}\partial_{y_2}\widehat \varphi_1,\\
\\
\displaystyle 
\partial_{z_1}\widetilde\varphi_{\varepsilon,2}(z)=\partial_{z_1}\widehat \varphi_2+ \varepsilon^{-1}\partial_{y_1}\widehat \varphi_2,\quad h_\varepsilon^{-1}\partial_{z_2}\widetilde\varphi_{\varepsilon,1}(z)= \varepsilon^{-1}\partial_{y_2}\widehat \varphi_2,\\
\end{array}$$
which can be written as follows
$$\begin{array}{l}\displaystyle
\mathbb{D}_{h_\varepsilon}[\widetilde \varphi_{\varepsilon}(z)]=\mathbb{D}_{h_\varepsilon}[\widehat \varphi]+\varepsilon^{-1}\mathbb{D}_{y}[\widehat \varphi],\quad 
{\rm div}_{h_\varepsilon}(\widetilde \varphi_\varepsilon(z))=\partial_{z_1}\widehat \varphi_1 +\varepsilon^{-1}{\rm div}_y(\widehat \varphi ).
\end{array}
$$
Then, we have that (\ref{form_var_general_1}) reads as follows
\begin{equation}\label{form_var_hat_u0}
\begin{array}{l}
\displaystyle \nu\int_{\widetilde\Omega_\ep}S_r(\mathbb{D}_{h_\ep}[\widetilde v_\ep]): \mathbb{D}_{h_\ep}[\widehat \varphi]\,dz  +{\nu\ep^{-1}}\int_{\widetilde\Omega_\ep}S_r(\mathbb{D}_{h_\ep}[\widetilde v_\ep]): \mathbb{D}_y[\widehat \varphi]\,dz  \\
\\
\displaystyle   - \int_{\Omega}\widetilde P_\ep^1\,\partial_{z_1}\widehat \varphi_1\,dz-\ep^{-1}\int_{\Omega }\widetilde P_\ep^1\,{\rm div}_{y}(\widehat \varphi)\,dz=
\int_{\widetilde\Omega_\ep}f_1\, \widehat \varphi_1\,dz\,.
\end{array} 
\end{equation}
Applying H${\rm \ddot{o}}$lder's inequality and taking into account estimates (\ref{estimates_critical_case}) and $\ep\ll h_\ep$ given in (\ref{parameters}), we get
$$\begin{array}{l}
\displaystyle 
\left|  \nu\int_{\widetilde\Omega_\ep}S_r(\mathbb{D}_{h_\ep}[\widetilde v_\ep]): \mathbb{D}_{h_\ep}[\widehat \varphi]\,dz\right|\leq   C\ep\|D_{h_\ep}\widehat \varphi\|_{L^{r}(\widetilde\Omega_\ep)^{2\times 2}}\leq C\ep \|D_{z}\widehat \varphi\|_{L^{r}(\widetilde\Omega_\ep)^{2\times 2}}\leq C\ep h^{-1}_\ep\to 0,
\end{array}$$
and taking into account that $\ep^{-r}S_r(\mathbb{D}_y[\widehat v_\ep])=S_r(\ep^{-{r\over r-1}}\mathbb{D}_y[\widehat v_\ep])$, by the change of variables given in Remark \ref{remarkCV}, we obtain 
\begin{equation}\label{form_var_hat_u}
\begin{array}{l}
\displaystyle\nu\int_{\Omega\times Y_f}S(\ep^{-{r\over r-1}}\mathbb{D}_{y}[\widehat v_\ep]): \mathbb{D}_y\widehat \varphi\,dzdy   - \int_{\Omega\times Y}\widehat P_\ep^1\,\partial_{z_1}\widehat \varphi_1\,dzdy\\
\\
\displaystyle-\ep^{-1}\int_{\Omega\times Y}\widehat P_\ep^1\,{\rm div}_{y}(\widehat \varphi)\,dx'dz_3dy=
\int_{\Omega\times Y_f}f_1\, \widehat \varphi_1\,dzdy+O_\ep\,.
\end{array} 
\end{equation}

{\it Step 2. Passing to the limit. } Now, we want to prove that the pair of limit functions $(\widehat v, p^1)$   given in Lemmas \ref{lemma_compactness} and \ref{lemma_conv_pressure},  satisfies the following two pressure limit system 
\begin{equation}\label{system_1_2_hat}\left\{\begin{array}{rl}
-\nu\,{\rm div}_y(|\mathbb{D}_y[\widehat v]|^{r-2}\mathbb{D}_y[\widehat v])+ \nabla_y\widehat \pi= \left(f_1(z_1) -\partial_{z_1} p^1(z_1)\right){\rm e}_1 & \hbox{in }\quad   Y_f, \\
\\\displaystyle
{\rm div}_y(\widehat v) =0& \hbox{in }\quad   Y_f, 
\\
\\\displaystyle
 \widehat v=0 \hbox{ on }  T & \hbox{ for a.e. } z\in \Omega,
\\ \displaystyle {\partial}_{z_1}\left(\int_0^1\int_{Y_f} \widehat v_1(z,y)\, dydz_2\right)=0
&\hbox{ in }\omega,\\
\displaystyle \left(\int_0^1\int_{Y_f} \widehat v_1(z,y)\,dydz_2 \right) n=0 & \hbox{ on }\partial\omega,
 \\
 (\widehat v, \widehat \pi) \hbox{ is }Y-{\rm periodic},
\end{array}\right.
\end{equation}
which has a unique solution $(\widehat v,  p^1, \widehat \pi)\in L^r(\Omega;W^{1,r}_{\rm per}(Y_f)^2)\times (L^{r'}_0(\omega)\cap W^{1,r'}(\omega))\times L^{r'}(\Omega;L^{r'}_{0,{\rm per}}(Y_f))$.

To do this, we consider $\widehat w(z,y)=(\widehat w_1(z,y),\widehat w_2(z_1,y))\in \mathcal{D}(\Omega;W_{{\rm per}}^{1,r}(Y)^3)$,  such that $\widehat w=0$ in $\Omega\times T$  and ${\rm div}_y(\widehat w)=0$ in $\Omega\times Y_f$. Thus, we consider the following test function in  (\ref{form_var_hat_u}):
$$\widehat \varphi_\ep(z,y)=\widehat w(z,y)-\ep^{-{r\over r-1}}\widehat v_\ep.$$
So we have
$$
\begin{array}{l}
\displaystyle \nu\int_{\Omega\times Y_f}S_r(\ep^{-{r\over r-1}}\mathbb{D}_{y}[\widehat v_\ep]): \mathbb{D}_y[\widehat \varphi_\ep]\,dzdy-\int_{\Omega\times Y}\widehat P^1_\ep\,\partial_{z_1}\widehat \varphi_{\ep,1}\,dzdy=
\int_{\Omega\times Y_f}f_1\, \widehat \varphi_{\ep,1}\,dzdy+O_\ep\,.
\end{array} 
$$
 which is equivalent to
$$
\begin{array}{l}
\displaystyle \nu\int_{\Omega\times Y_f}\left(S_r(\mathbb{D}_{y}[\widehat w])-S_r(\ep^{-{r\over r-1}}\mathbb{D}_{y}[\widehat v_\ep])\right): \mathbb{D}_y[\widehat \varphi_\ep]\,dzdy-\nu\int_{\Omega\times Y_f} S_r(\mathbb{D}_{y}[\widehat w]): \mathbb{D}_y[\widehat \varphi_\ep]\,dzdy\\
\\
\displaystyle+\int_{\Omega\times Y}\widehat P^1_\ep\,\partial_{z_1}\widehat \varphi_{\ep,1}\,dzdy=-
\int_{\Omega\times Y_f}f_1\, \widehat \varphi_{\ep,1}\,dzdy+O_\ep\,. 
\end{array} 
$$
Since $S_r$ is monotone, i.e.
$$\left(S_r(\mathbb{D}_{y}[\widehat w])-S_r(\ep^{-{r\over r-1}}\mathbb{D}_{y}[\widehat v_\ep])\right): (\mathbb{D}_{y}[\widehat w(z,y)]-\ep^{-{r\over r-1}}\mathbb{D}_{y}[\widehat v_\ep]) \geq 0,$$
we can deduce
$$
\begin{array}{l}
\displaystyle \nu\int_{\Omega\times Y_f}S_r(\mathbb{D}_{y}[\widehat w]): \mathbb{D}_y[\widehat \varphi_\ep]\,dzdy-\int_{\Omega\times Y}\widehat P^1_\ep\,\partial_{z_1}\widehat \varphi_{\ep,1}\,dzdy\geq 
\int_{\Omega\times Y_f}f_1\, \widehat \varphi_{\ep,1}\,dzdy+O_\ep\,. 
\end{array} 
$$
Passing to the limit by using convergences (\ref{conv_vel_gorro}) and (\ref{conv_pressure_gorro}), we obtain 
$$
\begin{array}{l}
\displaystyle\nu\int_{\Omega\times Y_f}S_r(\mathbb{D}_{y}[\widehat w]): \mathbb{D}_y[\widehat w-\widehat v]\,dzdy-\int_{\Omega\times Y}p^1 \,\partial_{z_1}(\widehat w_{1}-\widehat v_1)\,dzdy \geq 
\int_{\Omega\times Y_f}f_1\, (\widehat w_1-\widehat v_1)\,dzdy\,. 
\end{array} 
$$
From Minty's Lemma [31, Chapter 3, Lemma 1.2], then previous
inequality is equivalent to the following variational formulation
\begin{equation}\label{limit_form_var}
\begin{array}{l}\displaystyle \nu \int_{\Omega\times Y_f}S_r(\mathbb{D}_{y}[\widehat v]) :D_{y}\widehat w\,dzdy-\int_{\Omega\times Y}p^1\,\partial_{z_1}\widehat w_{1}\,dzdy=
\int_{\Omega\times Y_f}f_1\,\widehat w_1\,dzdy\,.
\end{array}
\end{equation}
By density, this equality holds for every function in the Banach space $\mathcal{W}$  defined by 
$$
\mathcal{W}=\left\{\begin{array}{l}
\widehat w(z,y)\in L^r(\Omega;W^{1,r}_{\rm per}(Y)^2) \ : \ {\rm div}_y\widehat w(z,y)=0\hbox{ in }\Omega\times Y_f,\quad \widehat w(z,y)=0\hbox{ in }\Omega\times T
\end{array}\right\}\,.
$$

Reasoning as in  \cite{Bourgeat1}, by integration by parts, the variational formulation (\ref{limit_form_var}) is equivalent to the system (\ref{system_1_2_hat}), where $\widehat \pi$ arises as a Lagrange multiplier of the incompressibility
constraint ${\rm div}_y(\widehat w)$ in $\Omega\times Y_f$. From  \cite[Theorem 8]{Bourgeat1} (whose proof is similar to the proof of  \cite[Theorem 2]{Bourgeat1}), it holds the uniqueness of solution $(\widehat v,  p^1, \widehat \pi)\in L^r(\Omega;W^{1,r}_{\rm per}(Y_f)^2)\times (L^{r'}_0(\omega)\cap W^{1,r'}(\omega))\times L^{r'}(\Omega;L^{r'}_{0,{\rm per}}(Y_f))$ and so, the whole sequence converges.\\

{\it Step 3. Average velocity.} We will deduce the expression for velocity $V_{av}$ given in (\ref{VelProm}). To do this, let us define the local problems which are useful to eliminate the variable $y$ of the homogenized problem  (\ref{system_1_2_hat}).

Let us separate variables $y$ and $z$ in the homogenized problem (\ref{system_1_2_hat}) satisfied by $(\widehat v, p^1 , \widehat \pi)$. It is easy to check that $(\widehat v, \widehat \pi)$ can be written as follow
$$\begin{array}{l}
\displaystyle 
\widehat v(z,y)={1\over \nu^{r'-1}}\Big|f_1(z_1)-{d\over d{z_1}}p^1(z_1)\Big|^{r'-2}\Big(f_1(z_1)-{d\over d{z_1}}p^1(z_1)\Big) \widehat w(y),
\\
\\
\displaystyle \widehat \pi(z,y)=\Big(f_1(z_1)-{d\over d{z_1}}p^1(z_1)\Big)\widehat q(y),
\end{array}$$
where $(\widehat w, \widehat q)$ is the unique solution of problem (\ref{Local_problems}).\\

Finally, taking into account the definition of  $V_{av}$ and relation (\ref{relation_u_ugorro}), we have 
$$\begin{array}{l}\displaystyle 
V_{av,1}(z_1)={1\over \nu^{r'-1}}\Big|f_1(z_1)-{d\over d{z_1}}p^1(z_1)\Big|^{r'-2}\Big(f_1(z_1)-{d\over d{z_1}}p^1(z_1)\Big) \int_{Y_f}\widehat w(y)\cdot {\rm e}_1dy\\
\\
V_{av,2}(z_1)\equiv 0,
\end{array}$$
and from relation (\ref{relation_u_ugorro}), we have (\ref{property_pressure}). Taking into account the definition of $\mu$ given in (\ref{nueff}), gives expression (\ref{VelProm}). 

\end{proof}
\section{Critical case: problem in the thin film}\label{sec:thinfilm}
As in the previous section, we assume the relation between the parameters in the critical regime (\ref{CritCase}). From this and  $\ep\ll \eta_\ep$ given in (\ref{parameters}), we deduce that  the last two terms in the estimate of $\mathcal{\widetilde U}_\ep$ given in (\ref{estim_u_I_dil}) satisfy
$$\begin{array}{l}
\displaystyle \ep^{1\over r-1}h_\ep^{1\over r}\eta_\ep^{r-1\over r}\approx \ep^{1\over r-1}\eta_\ep^{r-1\over r}\eta_\ep^{2r-1\over r(r-1)} \ep^{-{1\over r-1}}\ll \eta_\ep^{{1\over r-1}+{r-1\over r}+{2r-1\over r(r-1)}-{1\over r-1}}= \eta_\ep^{r\over r-1},\\
\end{array}$$
and 
$$\begin{array}{l}
\displaystyle \ep^{r+1\over 2(r-1)}h_\ep^{1\over r}\eta_\ep^{r-2\over 2r}\approx  \ep^{r+1\over 2(r-1)}\eta_\ep^{r-2\over 2r}\eta_\ep^{2r-1\over r(r-1)} \ep^{-{1\over r-1}}\ll \eta_\ep^{{r+1\over 2(r-1)}+{r-2\over 2r}+{2r-1\over r(r-1)}-{1\over r-1}}= \eta_\ep^{r\over r-1}.\\
\end{array}$$
Then, estimates (\ref{estim_u_I_dil})--(\ref{estim_Du_I_dil}) for $\mathcal{\widetilde U}_\ep$ and (\ref{esti_P_thin_film}) for $\widetilde P_\ep^2$ in $\widetilde I_1$ read as follows
\begin{equation}\label{estim_u_I_dil_crit}
\|\mathcal{\widetilde U}_\ep\|_{L^r(\widetilde I_1)^{2}}\leq C  \eta_\ep^{r\over r-1},\quad \|\mathbb{D}_{\eta_\ep}[\mathcal{\widetilde U}_\ep]\|_{L^r(\widetilde I_1)^{2\times 2}}\leq C\eta_\ep^{1\over r-1}, \quad
\|D_{\eta_\ep}\mathcal{\widetilde U}_\ep\|_{L^r(\widetilde I_1)^{2\times 2}}\leq C\eta_\ep^{1\over r-1},
\end{equation}
\begin{equation}\label{esti_P_thin_film_crit}
\|\widetilde P_\ep^2\|_{L^{r'}(\widetilde I_1)}\leq C,\quad 
\|\nabla_{\eta_\ep}  \widetilde P_\ep^2\|_{W^{-1,r'}(\widetilde I_1)^2}\leq C.
\end{equation}

\subsection{Convergences of velocity and pressure}
Using estimates  (\ref{estim_u_I_dil_crit}) and (\ref{esti_P_thin_film_crit}) and compactness, we prove the following lemma. 
\begin{lemma}\label{conv_tilde_func}
For a subsequence of $\ep$ still denoted by $\ep$, there exist:
\begin{itemize}
\item[--] $\mathcal{U}\in W^{1,r}(-g(z_1),0;L^r(\omega)^2)$, with $\mathcal{U}_2\equiv 0$ and $\mathcal{U}_1=0$ on $\Sigma\cup \Gamma_g$,  such that
\begin{equation}\label{conv_UP2_vel}
\eta_\ep^{-{r\over r-1}}\mathcal{\widetilde U}_\ep\rightharpoonup \mathcal{U}\quad\hbox{in }W^{1,r}(-g(z_1),0;L^r(\omega)^2),
\end{equation}
\item[--] $p^2\in L^{r'}_0(\widetilde I_1)\cap W^{1,r'}(\omega)$ independent of $z_2$,  such that
\begin{equation}\label{conv_UP2_press}
\widetilde P_\ep^2\rightharpoonup p^2\quad\hbox{in }L^{r'}(\widetilde I_1).
\end{equation}
\end{itemize}
\end{lemma}
\begin{proof}
From estimates (\ref{estim_u_I_dil_crit}), there exist $\mathcal{U}\in W^{1,r}(-g(z_1),0;L^r(\omega)^2)$ such that convergence (\ref{conv_UP2_vel}) holds.

Let $\widetilde\varphi\in C^\infty_0(\widetilde I_1)$, then 
$$\begin{array}{l}
\displaystyle \eta_\ep^{-{1\over r-1}}\int_{\widetilde I_1}\left(\partial_{z_1}\mathcal{\widetilde U}_{\ep,1}+\eta_\ep^{-1}\partial_{z_2}\mathcal{\widetilde U}_{\ep,2}\right)\widetilde\varphi\,dz
=-\eta_\ep^{-{1\over r-1}}\int_{\widetilde I_1}\mathcal{\widetilde U}_{\ep,1}\partial_{z_1}\widetilde\varphi\,dz-\eta_\ep^{-{r\over r-1}}\int_{\widetilde I_1}\mathcal{\widetilde U}_{\ep,2}\partial_{z_2}\widetilde\varphi\,dz.$$
\end{array}$$
Taking the limit $\ep\to 0$, we get
$$\int_{\widetilde I_1}\mathcal{U}_{2}\partial_{z_2}\widetilde\varphi\,dz=0,$$
so that $\mathcal{ U}_2$ does not depend on $z_2$.

Since $\mathcal{ U}$, $\partial_{z_2}\mathcal{U}\in L^r(\widetilde I_1)^2$, the traces $\mathcal{ U}(z_1,-g(z_1))$, $\mathcal{\widetilde U}(z_1,0)$  are well defined in $L^r(\omega)^2$. The proof of $\mathcal{ U}(z_1,-g(z_1))=0$ straightforward from the boundary condition $\mathcal{\widetilde U}_\ep(z_1, -g(z_1))=0$. Next, we prove that $\mathcal{U}(z_1,0)=0$. Proceeding similarly to the proof of Lemma \ref{Lem:Iep} (but now choosing a point $\beta_{z_1}\in T_\ep$ which is close to the point $\alpha_{z_1}\in \Sigma$), then we have
$$|\mathcal{\widetilde U_\ep}(\alpha_{z_1})|=|\widetilde u_\ep(\alpha_{z_1})-\widetilde u_\ep(\beta_{z_1})|=|D\widetilde u_\ep(\xi)(\alpha_{z_1}-\beta_{z_1})|\leq \ep |D_{h_\ep}\widetilde u_\ep|.$$
Since $\widetilde u_\ep(\beta_{z_1})=0$ because $\beta_{z_1}\in T_\ep$, we have
$$\|\mathcal{\widetilde  U}_\ep\|_{L^r(\Sigma)^2}\leq C\ep\|D_{h_\ep}\widetilde u_\ep\|_{L^r(\widetilde \Omega_\ep)^{2\times 2}}.$$
Then, taking into account estimate (\ref{estimates_critical_case}) in the critical case, we have
$$\|\eta_\ep^{-{r\over r-1}}\mathcal{\widetilde U}_\ep\|_{L^r(\Sigma)^2}\leq C\ep\eta_\ep^{-{r\over r-1}} \ep^{1\over r-1}=C(\ep/\eta_\ep)^{r\over r-1},$$
which tends to zero as $\ep$ to zero, because   $\ep\ll \eta_\ep$. This implies that $\mathcal{\widetilde U}(z_1,0)=0$.   Consequently, $\mathcal{\widetilde U}_2\equiv 0$, which finishes the proof of convergence (\ref{conv_UP2_vel}).
\\

Next, estimate (\ref{esti_P_thin_film_crit})$_1$ implies the existence of $p^2\in L^{r'}(\widetilde I_1)$ such that convergence (\ref{conv_UP2_press}) holds. Similarly to the proof of Lemma \ref{lemma_conv_pressure},  from estimate  (\ref{esti_P_thin_film_crit})$_2$, we deduce that $p_2$ does not depend on $z_2$. Since $\widetilde P^2$ has null mean value in $\widetilde I_1$, then $p^2$ also has null mean value in $\widetilde I_1$.

Now, by considering $\widetilde \varphi\in \mathcal{D}(\omega)$  as test function in the divergence condition ${\rm div}_{\eta_\ep}(\mathcal{\widetilde U}_\ep)=0$ in $\widetilde I_1$, we get
$$\int_{\widetilde I_1}\left({\partial}_{z_1}\mathcal{\widetilde U}_{\ep,1}\, \widetilde \varphi+\eta_\ep^{-1}\partial_{z_2}\mathcal{\widetilde U}_{\ep,2}\widetilde  \varphi\right)\,dz=0,$$
which, after integration by parts and multiplication by $\eta_\ep^{-{r\over r-1}}$, gives
$$\int_{\widetilde I_1} \eta_\ep^{-{r\over r-1}}\mathcal{\widetilde U}_{\ep,1}\, \partial_{z_1}\widetilde \varphi\,dz=0.$$
Passing to the limit by using convergence (\ref{conv_UP2_vel}), we deduce
$$\int_{\widetilde I_1}\mathcal{U}_1\, \partial_{z_1}\widetilde \varphi\,dz=0,$$
and, since $\varphi$ does not depend on $z_2$, we obtain the following divergence condition 
\begin{equation}\label{divxproperty_thin}
{\partial}_{z_1} \left(\int_{-g(z_1)}^0\mathcal{U}_1(z)\,dz_2\right)=0\quad \hbox{in }\omega,\quad   \left(\int_{-g(z_1)}^0\mathcal{U}_1(z)\,dz_2\right)\, n=0\quad \hbox{on }\partial\omega.
\end{equation}
By using convergences (\ref{conv_UP2_vel}) and (\ref{conv_UP2_press}), we refer to \cite[Propositions 3.1 and 3.2]{MikelicTapiero} in order to identify the effective system 
\begin{equation}\label{system_Utilde}
\left\{\begin{array}{rl}
\displaystyle -\partial_{z_2}\left(|\partial_{z_2}\mathcal{U}_1|^{r-2}\partial_{z_2}\mathcal{U}_1\right)={2^{r\over 2}\over \nu}\left(f_1(z_1)-\partial_{z_1}p^2(z_1)\right)& \hbox{in }\widetilde I_1,\\
\\
\displaystyle \partial_{z_1}\left(\int_{-g(z_1)}^0\mathcal{U}_1(z)\,dz_2\right)=0&\hbox{in }\omega,\\
\displaystyle
\left(\int_{-g(z_1)}^0\mathcal{U}_1(z)\,dz_2\right)\,n=0&\hbox{on }\partial\omega,\\
\\
\mathcal{U}_1=0&\hbox{ on }\Sigma\cup \Gamma_g.
\end{array}\right.
\displaystyle
\end{equation}
Furthermore, we have that $p^2\in W^{1,r'}(\omega)\cap L^{r'}_0(\widetilde I_1)$ according to  \cite[Proposition 3.3]{MikelicTapiero}. This concludes the proof.
\end{proof}
\subsection{Average velocity in the thin film}
\begin{theorem} Consider  $(\mathcal{U}, p^2)$   given in Lemma \ref{conv_tilde_func}. Then, we have that  the average velocity $$\mathcal{V}_{av}(z_1)={1\over g(z_1)}\int_{-g(z_1)}^0\mathcal{U}(z)\,dz_2,$$ is given by
\begin{equation}\label{Vexpav}
\left\{
\begin{array}{ll} 
\displaystyle \mathcal{V}_{av,1}(z_1)= {g(z_1)^{r'}\over 2^{r'\over 2}(r'+1)\nu^{r'-1}}\left|f_1(z_1)-{d\over dz_1}p^2(z_1)\right|^{r'-2}\left(f_1(z_1)-{d\over dz_1}p^2(z_1)\right)& \hbox{in }\omega,\\
\\
\displaystyle
\mathcal{V}_{av,2}\equiv 0.
\end{array}\right.
\end{equation}
\end{theorem}
\begin{proof} 
Since $\mathcal{U}_2\equiv 0$, it is only necessary to obtain the  expression of $\mathcal{U}_1$, which satisfies problem (\ref{system_Utilde})$_{1,3}$. We remark that this problem is formally an ordinary differential equation in the variable $z_2$, with parameter $z_1\in \omega$. The resolution of (\ref{system_Utilde}) is similar to \cite[Proposition 3.4]{MikelicTapiero} (see also \cite[Lemma 6.3]{Fabricius}), so we omit it. 
 
\end{proof}

 \section{Critical case:  a generalized Reynolds limit equation}\label{sec:mainresult}
 The conclusion of the previous two sections is that for any sequence of solutions $(\widetilde v_\ep, \widetilde P^{1}_\ep)$ and $(\mathcal{\widetilde U}_\ep, \widetilde P^2_\ep)$ and letting $\ep\to 0$, we can extract subsequences, still denoted by the same symbol, and find functions $(v, p^1)\in W^{1,r}(0, 1;L^r(\omega)^2)\times W^{1,r'}(\omega)$ and $(\mathcal{U}, p^2)\in W^{1,r}(-g(z_1),0;L^r(\omega)^2)\times W^{1,r'}(\omega)$ such that
 \begin{equation}\label{convergences_all}
 \begin{array}{ll}
 \displaystyle
\ep^{-{r\over r-1}}\widetilde V_\ep \rightharpoonup v=(v_1,0)\quad\hbox{in }W^{1,r}(0,1;L^r(\omega)^2),&\displaystyle \widetilde P_\ep^1\to p^1\quad\hbox{ in }L^{r'}(\Omega),\\
 \\
 \displaystyle
 \eta_\ep^{-{r\over r-1}}\mathcal{\widetilde U}_\ep\rightharpoonup \mathcal{U}=(\mathcal{U}_1,0)\quad\hbox{in }W^{1,r}(-g(z_1),0;L^r(\omega)^2),&\displaystyle \widetilde P_\ep^2\rightharpoonup p^2\quad\hbox{in }L^{r'}(\widetilde I_1).
 \end{array}
 \end{equation}
Moreover,   functions $(V_{av}, p^1), (\mathcal{V}_{av},p^2)$, with $V_{av}=\int_0^1v(z)\,dz_2$ and $\mathcal{V}_{av}=g(z_1)^{-1}\int_{-g(z_1)}^0\mathcal{U}(z)\,dz_2$, necessarily satisfy the following equations in $\omega$ 
 \begin{equation}\label{equations_all}
 \begin{array}{l}
 \displaystyle
 V_{av,1}(z_1)={\mu\over \nu^{r'-1}}\Big|f_1(z_1)-{d\over d{z_1}}p^1(z_1)\Big|^{r'-2}\Big(f_1(z_1)-{d\over d{z_1}}p^1(z_1)\Big),\quad V_{av,2}\equiv 0,\\
 \\
\displaystyle \mathcal{V}_{av,1}(z_1)= {g(z_1)^{r'}\over 2^{r'\over 2}(r'+1)\nu^{r'-1}}\left|f_1(z_1)-{d\over dz_1}p^2(z_1)\right|^{r'-2}\left(f_1(z_1)-{d\over dz_1}p^2(z_1)\right),\quad \mathcal{V}_{av,2}\equiv 0,
\end{array}
\end{equation}
 with $\mu$ defined in (\ref{nueff}).

 Next, we find the connection between the functions $p^1$ and $p^2$, i.e.  the
coupling effects between the solution in the porous  and free media.
 
 \begin{lemma} We assume that the parameters $\ep, \eta_\ep$ and $h_\ep$ satisfy (\ref{parameters}) and (\ref{CritCase}). Let   $p^1\in W^{1,r'}(\omega)$ and $p^2\in W^{1,r'}(\omega)$ be such that  (\ref{convergences_all}) and   (\ref{equations_all}) hold. Then, we have
 \begin{equation}\label{generalizedReynolds}
 \begin{array}{l}
 \displaystyle 
{\mu\over \nu^{r'-1}} \int_{\omega} \Big|f_1(z_1)-{d\over d{z_1}}p^1(z_1)\Big|^{r'-2}\Big(f_1(z_1)-{d\over d{z_1}}p^1(z_1)\Big)
 {d\over dz_1}\psi\,dz_1\\
 \\
 \displaystyle 
 +{1\over \lambda 2^{r'\over 2}(r'+1)\nu^{r'-1}}\int_{\omega}{g(z_1)^{r'}}\left|f_1(z_1)-{d\over dz_1}p^2(z_1)\right|^{r'-2}\left(f_1(z_1)-{d\over dz_1}p^2(z_1)\right){d\over dz_1}\psi \,dz_1=0,
 \end{array}
 \end{equation}
 for every $\psi\in W^{1,r'}(\omega)$.
  \end{lemma}
 \begin{proof}
 Choosing  $\psi\in W^{1,r'}(\omega)$ as test function in (\ref{form_var_1_tilde})$_2$, putting $\widetilde u_\ep=0$ in the solid part  (recalling that $\widetilde V_\ep$ is the extension by zero of $\widetilde v_\ep$ to the whole $\Omega$) and integrating by parts, we get
  \begin{equation}\label{coupling1}
 \int_{\Omega}h_\ep  \widetilde V_{\ep,1} \partial_{z_1}\psi(z_1)\,dz+\int_{\widetilde I_1}\eta_\ep\,  \mathcal{\widetilde U}_{\ep,1} \partial_{z_1}\psi(z_1)\,dz=0.
  \end{equation}
Multiplying (\ref{coupling1}) by $\ep^{-{r\over r-1}}h_\ep^{-1}$, we have
 $$\int_{\Omega}\ep^{-{r\over r-1}} \widetilde V_{\ep,1} \partial_{z_1}  \psi(z_1)\,dz+\int_{\widetilde I_1}\ep^{-{r\over r-1}}h_\ep^{-1}\eta_\ep  \mathcal{\widetilde U}_{\ep,1} \partial_{z_1} \psi(z_1)\,dz=0,$$
which can be rewritten as follows 
  $$\int_{\Omega}\ep^{-{r\over r-1}} \widetilde V_{\ep,1} \partial_{z_1}  \psi(z_1)\,dz+\ep^{-{r\over r-1}}\eta_\ep^{2r-1\over r-1}h_\ep^{-1}\int_{\widetilde I_1}\eta_\ep^{-{r\over r-1}}\mathcal{\widetilde U}_{\ep,1} \partial_{z_1}  \psi(z_1)\,dz=0.$$
From convergences (\ref{convergences_all}) and relation (\ref{CritCase}),   passing to the limit as $\ep\to 0$, we deduce
 $$\int_{ \Omega }v_1 \partial_{z_1}  \psi(z_1)\,dz+\lambda^{-1}\int_{\widetilde I_1} \mathcal{U}_1 \partial_{z_1} \psi(z_1)\,dz=0.$$
 Then, since $\psi$ does not depend on $z_2$, we have 
$$ \int_{\omega} V_{av,1}\partial_{z_1}\psi(z_1)\,dz_1+\lambda^{-1}\int_{\omega}g(z_1)\mathcal{ V}_{av, 1} \partial_{z_1}\psi(z_1)\,dz_1=0.$$
Taking into account (\ref{equations_all}), this is the equation (\ref{generalizedReynolds}).
 \end{proof}

 In the following result, we are going to prove the relation between the pressures $p^1$ and $p^2$, i.e. the continuity of the pressure in $\Sigma$.

 \begin{lemma}We assume that the parameters $\ep, \eta_\ep$ and $h_\ep$ satisfy (\ref{parameters}) and (\ref{CritCase}). Let $p^1$ and $p^2$ be the limit pressures from expression (\ref{convergences_all}). Then, there exists $c^\star\in \mathbb{R}$ such that
\begin{equation}\label{relation_pressures}
p^1=p^2+c^*.
\end{equation}
 \end{lemma}
 \begin{proof}  For any $\widetilde \varphi\in \mathcal{D}(\omega)$, we define $(\widetilde \phi, \widetilde \psi)\in W^{1,r}(\Omega)\times W^{1,r}(\widetilde I_1)$ such that
 $$\widetilde \phi=0\quad\hbox{on } \widetilde \Omega\setminus \Sigma,\quad \widetilde \psi=0\quad\hbox{on }\partial \widetilde I_1\setminus \Sigma,\quad \widetilde \phi=\widetilde \psi=\widetilde \varphi\quad \hbox{ on }\Sigma.$$
Let us define the following global test function in $\widetilde D_\ep$ given by
 \begin{equation}\label{wtilde}\widetilde w_\ep(z)=\left\{\begin{array}{ll}
\widetilde \phi (z)(\mathcal{\widetilde R}^\ep_r{\rm e_2})(z)&\hbox{in }\widetilde \Omega_\ep,\\
 \\
 \widetilde \psi(z){\rm e}_2&\hbox{in }\widetilde I_1,
 \end{array}\right.
 \end{equation}
 where $\mathcal{\widetilde R}^\ep_r$ is the restriction operator defined in Lemma \ref{restriction_operator2}.
 We observe that $ \mathcal{R}^\ep_r{\rm e}_2$ tends to its $Y$-average $\int_Y   (\mathcal{R}_r{\rm e}_2)(y)\,dy$ (where the restriction operator $\mathcal{R}_r$ is defined in the proof of Lemma \ref{restriction_operator}), and $\mathcal{\widetilde R}^\ep_r({\rm e}_2)_2$ tends to $1$ in $L^{r}(\Omega)$. 
\\
  
\noindent Now, we take $\widetilde w_\ep$ as test function  in the system (\ref{form_var_1_tilde}), and we obtain
 \begin{equation}\label{form_var_1_tilde_final}
 \begin{array}{l}
\displaystyle \nu\int_{\widetilde \Omega_\ep}h_\ep S_r(\mathbb{D}_{h_\ep}[\widetilde v_\ep]):\mathbb{D}_{h_\ep}[\widetilde w_\ep]\,dz+\nu\int_{\widetilde I_1}\eta_\ep S_r(|\mathbb{D}_{\eta_\ep}[\mathcal{\widetilde U}_\ep]):\mathbb{D}_{\eta_\ep}[\widetilde w_\ep]\,dz\\
\\
\displaystyle  -\int_{\widetilde \Omega_\ep}h_\ep\widetilde p_\ep\,{\rm div}_{h_\ep}(\widetilde w_\ep)\,dz-\int_{\widetilde I_1}\eta_\ep \widetilde p_\ep\,{\rm div}_{\eta_\ep}(\widetilde w_\ep)\,dz\\
\\
\displaystyle =\int_{\widetilde \Omega_\ep}h_\ep f\cdot\widetilde w_\ep\,dz+\int_{\widetilde I_1}\eta_\ep f\cdot\widetilde w_\ep\,dz.
\end{array} 
\end{equation}
From H${\rm \ddot{o}}$lder's inequality, estimates  (\ref{estim_restricted2}) and (\ref{estimates_critical_case}), and $\ep\ll h_\ep$, we deduce
 \begin{equation}\label{B1}
\begin{array}{l}\displaystyle 
\left|\displaystyle \nu\int_{\widetilde \Omega_\ep}h_\ep S_r(|\mathbb{D}_{h_\ep}[\widetilde v_\ep]):\mathbb{D}_{h_\ep}[\widetilde w_\ep]\,dz\right|\\
\\
\displaystyle =\left|\displaystyle \nu h_\ep\int_{\widetilde \Omega_\ep} S_r(|\mathbb{D}_{h_\ep}[\widetilde v_\ep]):(\nabla_{h_\ep}\widetilde\phi(z)\cdot(\mathcal{R}^\ep_r {\rm e}_2)(z))\,dz+\nu h_\ep \int_{\widetilde \Omega_\ep} S_r(|\mathbb{D}_{h_\ep}[\widetilde v_\ep]):\mathbb{D}_{h_\ep}[(\mathcal{\widetilde R}^\ep_r {\rm e}_2)(z)]\widetilde\phi(z)\,dz\right|
\\
\\
\displaystyle \leq C h_\ep\left(\|\mathbb{D}_{h_\ep}[\widetilde v_\ep]\|^{r-1}_{L^{r}(\widetilde\Omega_\ep)^{2\times 2}} \|\nabla_{h_\ep}\widetilde\phi\|_{L^r(\Omega)^2}+ \|\mathbb{D}_{h_\ep}[\widetilde v_\ep]\|^{r-1}_{L^{r}(\widetilde\Omega_\ep)^{2\times 2}} \|D_{h_\ep}\mathcal{\widetilde R}^\ep_r {\rm e}_2\|_{L^r(\Omega)^2}\right)\\
\\
\displaystyle \leq Ch_\ep(\ep h_\ep^{-1} +\ep \ep^{-1})\leq Ch_\ep,
\end{array}
\end{equation}
and by using estimates (\ref{estim_u_I_dil_crit}), we deduce
\begin{equation}\label{B2}\begin{array}{l}\displaystyle 
\left|\nu\int_{\widetilde I_1}\eta_\ep S_r(|\mathbb{D}_{\eta_\ep}[\mathcal{\widetilde U}_\ep]):\mathbb{D}_{\eta_\ep}[\widetilde w_\ep]\,dz\right|
=\displaystyle \left|\nu   \int_{\widetilde I_1} S_r(|\mathbb{D}_{\eta_\ep}[\mathcal{\widetilde U}_\ep]):\left(\eta_\ep\partial_{z_1}[\widetilde\psi(z){\rm e}_2]+ \partial_{z_2}[\widetilde\psi(z){\rm e}_2]\right)\,dz\right|\\
\\
\leq \displaystyle  C(\eta_\ep+1)\|\mathbb{D}_{\eta_\ep}[\mathcal{\widetilde U}_\ep]\|_{L^r(\widetilde I_1)^{2\times 2}}^{r-1}\leq C\eta_\ep.
\end{array}
\end{equation}
From the unfolding change of variables (\ref{CV}) and   ${\rm div}_y(\mathcal{\widetilde R}^\ep_r{\rm e_2})=0$ in $Y$, we deduce that
\begin{equation}\label{B2}\begin{array}{l}\displaystyle 
\int_{\widetilde \Omega_\ep}h_\ep\widetilde p_\ep\,{\rm div}_{h_\ep}(\widetilde w_\ep)\,dz
=\displaystyle \int_{\Omega}h_\ep \widetilde P^1_\ep\,{\rm div}_{h_\ep}(\widetilde \phi(z)(\mathcal{\widetilde R}^\ep_r{\rm e_2})(z))\,dz\\
\\
  =\displaystyle \int_{\Omega}h_\ep \widetilde P^1_\ep\,\nabla_{h_\ep}\widetilde \phi(z)\cdot (\mathcal{\widetilde R}^\ep_r{\rm e_2})(z)\,dz+ \int_{\Omega\times Y}h_\ep\ep^{-1} \widehat P^1_\ep\,\widetilde\phi(z)\,{\rm div}_{y}(\mathcal{\widetilde R}^\ep_r({\rm e_2}))\,dzdy\\
\\
=\displaystyle\displaystyle  \int_{\Omega}h_\ep \widetilde P^1_\ep\,\nabla_{h_\ep}\widetilde \phi(z)\cdot (\mathcal{\widetilde R}^\ep_r{\rm e_2})(z)\,dz,
\end{array}
\end{equation}
and from estimate (\ref{estim_P_hat}), we have
$$\left|\int_{\Omega}h_\ep \widehat P^1_\ep\,\nabla_{h_\ep}\widetilde \phi(z)\cdot (\mathcal{\widetilde R}^\ep_r{\rm e_2})(z)\,dz\right|\leq C.$$

\noindent  From the definition of $\widetilde P^2_\ep$ given in (\ref{P2}) and $\widetilde c_\ep$ given in (\ref{cep}),  
 we have that 
\begin{equation}\label{B4}\begin{array}{l}
 \displaystyle \eta_\ep  \int_{\widetilde I_1}\widetilde p_\ep\,{\rm div}_{\eta_\ep}(\widetilde w_\ep)\,dz
 =\displaystyle \eta_\ep  \int_{\widetilde I_1}(\widetilde p_\ep-\widetilde c_\ep)\,{\rm div}_{\eta_\ep}(\widetilde w_\ep)\,dz+\eta_\ep  \widetilde c_\ep\int_{\widetilde I_1}{\rm div}_{\eta_\ep}(\widetilde w_\ep)\,dz\\
\\ =\displaystyle    \int_{\widetilde I_1}\widetilde P^2_\ep\,\partial_{z_2}\widetilde \psi(z)\,dz+   \widetilde c_\ep\int_{\widetilde I_1}\partial_{z_2}\widetilde \psi(z)\,dz,
 \end{array}
 \end{equation}
and from estimate (\ref{esti_P_thin_film_crit}), we have
$$\left| \int_{\widetilde I_1}\widetilde P^2_\ep\,\partial_{z_2}\widetilde \psi(z)\,dz\right|\leq C.$$
\noindent  Taking into account that that $f$ is given by (\ref{fep}) and $\widetilde w_\ep$ given in (\ref{wtilde}), we deduce 
\begin{equation}\label{B5}
\begin{array}{l}
\displaystyle  \int_{\widetilde \Omega_\ep}h_\ep f\cdot\widetilde w_\ep\,dz=0,\quad\hbox{and}\quad
\displaystyle \int_{\widetilde I_1}\eta_\ep   f\cdot\widetilde w_\ep\,dz=0.
\end{array}
\end{equation}
\noindent From (\ref{B1})--(\ref{B5}) and the fact that $\eta_\ep\ll 1$ and $h_\ep\ll 1$, we deduce that  $|\widetilde c_\ep|\leq C$ and so there exists $c^*$ such that $\widetilde c_\ep$ tends to $c^*$. Moreover,  we deduce that (\ref{form_var_1_tilde_final}) reads as follow
\begin{equation}\label{form_var_1_tilde_final2}
 \begin{array}{l}
\displaystyle  \int_{\Omega}h_\ep \widetilde P^1_\ep\,\nabla_{h_\ep}\widetilde \phi(z)\cdot (\mathcal{\widetilde R}^\ep_r{\rm e_2})(z)\,dz+\int_{\widetilde I_1}\widetilde P^2_\ep\,\partial_{z_2}\widetilde \psi(z)\,dz+   \widetilde c_\ep\int_{\widetilde I_1}\partial_{z_2}\widetilde \psi(z)\,dz+O_\ep=0.
\end{array} 
\end{equation}
 Passing to the limit when $\ep\to 0$, from strong convergence of $\widetilde P^1_\ep$ given in (\ref{conv_pressure_sub}) and convergence $\mathcal{\widetilde R}^\ep_r{\rm e_2}$ to $1$ in the first term, convergence of $\widetilde P^2_\ep$ given in (\ref{conv_UP2_press}) in the second term and convergence of $\widetilde c_\ep$ to $c^\star$ in the third term, we get
$$
 \begin{array}{l}
\displaystyle \int_{\Omega} p^1(z_1)\partial_{z_2}\widetilde \phi(z)\,dz+\int_{\widetilde I_1}p^2(z_1)\,\partial_{z_2}\widetilde \psi(z)\,dz+   c^\star\int_{\widetilde I_1}\partial_{z_2}\widetilde \psi(z)\,dz=0.
\end{array} 
$$
Since $p^1$ and $p^2$ do not depend on $z_2$, this can be written as follows
$$
 \begin{array}{l}
\displaystyle \int_{\omega} p^1(z_1)\left(\int_0^1\partial_{z_2}\widetilde \phi(z)\,dz_2\right)dz_1+\int_{\omega}p^2(z_1)\left(\int_{-g(z_1)}^0\partial_{z_2}\widetilde \psi(z)dz_2\right)dz_1+ c^\star\int_{\omega}\int_{-g(z_1)}^0\partial_{z_2}\widetilde \psi(z)dz_2dz_1=0,
\end{array} 
$$
and integrating with respect to $z_2$, by taking into account that $\widetilde \phi(z_1, 1)=\widetilde\psi(z_1,-g(z_1))=0$, we get
$$-\int_\omega p^1(z_1)\widetilde \phi(z_1,0)\,dz_1+\int_{\omega}p^2(z_1)\,\widetilde \psi(z_1,0)\,dz_1+   c^\star\int_\omega\widetilde \psi(z_1,0)\,dz_1=0,$$
and taking into account that $\widetilde\phi=\widetilde\psi=\widetilde\varphi$ on $\Sigma$, then  we deduce
$$
-\int_\omega p^1(z_1)\widetilde \varphi(z_1)\,dz_1+\int_{\omega}(p^2(z_1)+c^\star)\,\widetilde \varphi(z_1)\,dz_1=0,
$$
for any $\widetilde\varphi\in \mathcal{D}(\omega)$, which implies that  equation (\ref{relation_pressures}) holds.

 \end{proof}
 
We have already proved the convergence of $\widetilde P^1_\ep$ to $p^1$ in $\Omega$ and $\widetilde P^2_\varepsilon$ to $p^2$ in $\widetilde I_1$. Let us define the following pressure in $D$ by
\begin{equation}\label{pstar}p^\star=\left\{
\begin{array}{ll}
p^1 &\hbox{in }\Omega,\\
p^2 +c^\star &\hbox{in }\widetilde I_1.
\end{array}\right.
\end{equation} Next, we give the main result of this paper. 
\begin{theorem}\label{main_thm}We assume  that the parameters $\ep, \eta_\ep$ and $h_\ep$ satisfy (\ref{parameters}) and (\ref{CritCase}). Then, the asymptotic pressure $p^\star$ defined in (\ref{pstar}) is the unique solution of the generalized Reynolds equation:
\\

\indent   Find $p^\star\in W^{1,r'}(\omega)\cap L^{r'}_0(\omega)$ such that 
\begin{equation}\label{generalizedReynolds_final}
 \begin{array}{l}
 \displaystyle  
\int_{\omega} \left({\mu\over \nu^{r'-1}} +{{g(z_1)^{r'}}\over \lambda 2^{r'\over 2}(r'+1)\nu^{r'-1}}\right)\left|f_1(z_1)-{d\over dz_1}p^\star(z_1)\right|^{r'-2}\left(f_1(z_1)-{d\over dz_1}p^\star(z_1)\right){d\over dz_1}\psi \,dz_1=0,
 \end{array}
\end{equation}
\indent for every  $\psi\in W^{1,r'}(\omega)$. 
\\

\noindent Moreover,   the average velocity field in the free media  is given by
$$
\left\{\begin{array}{l}
\displaystyle \mathcal{V}_{av,1}(z_1)= {g(z_1)^{r'}\over 2^{r'\over 2}(r'+1)\nu^{r'-1}}\left|f_1(z_1)-{d\over dz_1}p^\star(z_1)\right|^{r'-2}\left(f_1(z_1)-{d\over dz_1}p^\star(z_1)\right)\\
\\
\mathcal{V}_{av,2}\equiv 0\end{array}\right.
\quad \hbox{in }\omega,
$$ 
and the average velocity field in the porous media is given by
 $$\left\{\begin{array}{l}
 \displaystyle
 V_{av,1}(z_1)={\mu\over \nu^{r'-1}}\Big|f_1(z_1)-{d\over d{z_1}}p^\star(z_1)\Big|^{r'-2}\Big(f_1(z_1)-{d\over d{z_1}}p^\star(z_1)\Big)\\
 \\
 V_{av,2}\equiv 0,
\end{array}\right. \quad  \hbox{in }\omega,
 $$
with $\mu>0$  defined by 
\begin{equation}\label{mueff}
\mu=\int_{Y_f}|\mathbb{D}_y[\widehat w]|^r\,dy,
\end{equation}
where $(\widehat w, \widehat q)\in W^{1,r}_{\rm per}(Y)^2\times L^{r'}_{0, {\rm per}}(Y)$, $1<r<+\infty$, is the unique solution of the auxiliary problem 
\begin{equation}\label{boundary_layer}
\left\{\begin{array}{rl}
-{\rm div}_y\left(|\mathbb{D}_y[\widehat w]|^{r-2}\mathbb{D}_y[\widehat w]\right)+ \nabla_y \widehat q={\rm e}_1 & \hbox{ in }Y_f,\\
\\
{\rm div}_y \widehat w=0& \hbox{ in }Y_f.\\
\\
\widehat w=0& \hbox{ in }T.
\end{array}\right.
\end{equation}
 \end{theorem}
 \begin{proof}    All the results presented in the theorem are consequences of the previous results. In particular, from equation (\ref{generalizedReynolds}), (\ref{relation_pressures}) and (\ref{pstar}), we obtain the variational formulation for the limit pressure (\ref{generalizedReynolds_final}).

Let us prove that the uniqueness of  solution  up to an additive constant of (\ref{generalizedReynolds_final}). The proof relies on standard monotonicity arguments. Let us first introduce some notation and properties. Thanks to (\ref{Funcg}), we have  
\begin{equation}\label{estiminfG}
G(z_1):={\mu\over \nu^{r'-1}} +{{g(z_1)^{r'}}\over \lambda 2^{r'\over 2}(r'+1)\nu^{r'-1}}\geq {\mu\over \nu^{r'-1}} +{a^{r'}\over \lambda 2^{r'\over 2}(r'+1)\nu^{r'-1}}:=C_{\mu,\nu,r'}^{\lambda,a},
\end{equation}
with $C_{\mu,\nu,r'}^{\lambda,a}>0$ a  constant.  For $1<r'<+\infty$, we define the $r'$-Laplace operator $A_{r'}(\xi)=|\xi|^{r'-2}\xi$, $\forall \xi\in \mathbb{R}$, which is strongly monotone in the following sense (see for instance \cite{Baranger, Duvnjak}):
 \begin{itemize}
 \item If $1<r'< 2$, then there exists $\alpha_1>0$ such that
\begin{equation}\label{monsub}\begin{array}{l}\displaystyle
 \int_{\omega} \left(A_{r'}\left(u\right)-A_{r'}\left(v\right)\right)\left(u-v\right) \,dz_1
 \geq 
 \alpha_1{\|u-v\|^2_{L^{r'}(\omega)}\over(\|u\|_{L^{r'}(\omega)}+\|v\|_{L^{r'}(\omega)})^{2-r'}},
 \end{array}
 \end{equation}
 \item If  $r'\geq 2$, then there exists $\alpha_2>0$ such that
\begin{equation}\label{monsup}\begin{array}{l}\displaystyle
 \int_{\omega} \left(A_{r'}\left(u\right)-A_{r'}\left(v\right)\right) \left(u-v\right)  \,dz_1
 \geq 
 {\alpha_2\|u-v\|^{r'}_{L^{r'}(\omega)}}.
 \end{array}
 \end{equation}
 \end{itemize}Let us now suppose that  (\ref{generalizedReynolds_final}) has two solutions  $p,q\in W^{1,r'}(\omega)$. Then, subtracting the corresponding variational equations and taking $\psi(z_1)=q(z_1)-p(z_1)$ as test function, we get
 $$\begin{array}{l}
 \displaystyle  
\int_{\omega} G(z_1)\left\{A_{r'}\left(f_1(z_1)-{d\over dz_1}p(z_1)\right)-A_{r'}\left(f_1(z_1)-{d\over dz_1}q(z_1)\right)\right\}{d\over dz_1}(q(z_1)-p(z_1)) \,dz_1=0,
 \end{array}$$
 or, equivalently,
  $$\begin{array}{l}
 \displaystyle  
\int_{\omega} G(z_1)\left\{A_{r'}\left(f_1(z_1)-{d\over dz_1}p(z_1)\right)-A_{r'}\left(f_1(z_1)-{d\over dz_1}q(z_1)\right)\right\}\left\{\left(f_1-{d\over dz_1}p\right)-\left(f_1-{d\over dz_1}q\right) \right\}\,dz_1=0.
 \end{array}$$
Taking into account (\ref{estiminfG}), it holds
$$\begin{array}{l}
 \displaystyle C_{\mu,\nu,r'}^{\lambda,a}\int_{\omega} \left\{A_{r'}\left(f_1-{d\over dz_1}p\right)-A_{r'}\left(f_1-{d\over dz_1}q\right)\right\}\left\{\left(f_1-{d\over dz_1}p\right)-\left(f_1-{d\over dz_1}q\right)\right\} \,dz_1 
 \\
 \\
 \leq \displaystyle  
 \int_{\omega} G(z_1)\left\{A_{r'}\left(f_1-{d\over dz_1}p\right)-A_{r'}\left(f_1-{d\over dz_1}q\right)\right\}\left\{\left(f_1-{d\over dz_1}p\right)-\left(f_1-{d\over dz_1}q\right)\right\} \,dz_1=0,
 \end{array}$$
and then, 
\begin{equation}\label{final1}\begin{array}{l}
 \displaystyle C_{\mu,\nu,r'}^{\lambda,a}\int_{\omega} \left\{A_{r'}\left(f_1-{d\over dz_1}p\right)-A_{r'}\left(f_1-{d\over dz_1}q\right)\right\}\left\{\left(f_1-{d\over dz_1}p\right)-\left(f_1-{d\over dz_1}q\right)\right\} \,dz_1 \leq 0.
 \end{array}\end{equation}
 By respectively using the monotonicity properties (\ref{monsub}) for $1<r'<2$ and (\ref{monsup})  for $r'\geq 2$ applied to the left-hand side of (\ref{final1}), we deduce for  $r'\in (1,+\infty)$ that
 $$\begin{array}{l}
 \displaystyle \left\|{d\over dz_1}(p-q)\right\|_{L^{r'}(\omega)} = 0.
 \end{array}$$
 Thus, we get that $\partial_{z_1}(p(z_1)-q(z_1))=0$ in $\omega$ and so, the solution of (\ref{generalizedReynolds_final})  is only determined up to an additive constant. Taking into account (\ref{pstar}) and Lemma \ref{lemma_conv_pressure}, which says that $p^1\in L^{r'}_0(\omega)$,  a supplementary constraint has to be added to obtain $p^\star$, namely $p^\star\in W^{1,r'}(\omega)\cap L^{r'}_0(\omega)$ is the unique solution of (\ref{generalizedReynolds_final}).

 Finally, we just observe that multiplying equation (\ref{boundary_layer}) by $\widehat w$, integrating in $Y_f$ and taking into account that $\widehat w=0$ on $T$, we deduce that the permeability constant defined by (\ref{nueff}) also satisfies
 $$\mu=\int_{Y_f} \widehat w\cdot {\rm e}_1\,dy=\int_{Y_f}|\mathbb{D}_y[\widehat w]|^r\,dy,$$
which is (\ref{mueff}).
 \end{proof}
 
 \section{Conclusions}
{\bf Main result.} In this paper, we consider an incompressible viscous stationary 2D non-Newtonian fluid in a domain composed by two parts in contact: a periodic thin porous medium $\Omega_\ep$ with characteristic size of the pores $0<\ep\ll1$ and thickness of the domain $0<h_\ep\ll 1$, and a thin film $I_\ep$ with thickness $0<\eta_\ep\ll 1$, where $h_\ep$ and $\eta_\ep$ are devoted to zero when $\ep\to 0$.  The interface between $\Omega_\epsilon$ and $I_\epsilon$ is defined by $\Sigma=\omega\times \{x_2=0\}$. More precisely, we consider the case of a non-Newtonian fluid governed by the incompressible Stokes equations with power law viscosity of flow index $r\in (1, +\infty)$, and we prove that there exists a critical regime between these parameters given by 
$$
h_\ep\approx  \eta_\ep^{2r-1\over r-1} \ep^{-{r\over r-1}},\quad \hbox{i.e. }\quad {h_\ep\over  \eta_\ep^{2r-1\over r-1} \ep^{-{r\over r-1}}}\to \lambda\in (0,+\infty),
$$
where the pressure has the same order of magnitude in the porous medium and in the free film (with  a continuity relation of their limits through $\Sigma$) and is described by a modified Reynolds equation, coupling the effects of the thin porous medium (1D nonlinear Darcy problem with permeability $\mu>0$ given by (\ref{mueff})) and the thin film (1D nonlinear Reynolds problem), which given by
$$\left\{\begin{array}{ll}
\displaystyle -{d\over d z_1}\left[\left({\mu\over \nu^{r'-1}} +{{g(z_1)^{r'}}\over \lambda 2^{r'\over 2}(r'+1)\nu^{r'-1}}\right)\left|f_1(z_1)-{d\over dz_1}p^\star(z_1)\right|^{r'-2}\left(f_1(z_1)-{d\over dz_1}p^\star(z_1)\right)\right]=0 & \hbox{in }\omega,\\
\\
\displaystyle \left({\mu\over \nu^{r'-1}} +{{g(z_1)^{r'}}\over \lambda 2^{r'\over 2}(r'+1)\nu^{r'-1}}\right)\left|f_1(z_1)-{d\over dz_1}p^\star(z_1)\right|^{r'-2}\left(f_1(z_1)-{d\over dz_1}p^\star(z_1)\right)\, n=0 & \hbox{on }\partial\omega,
\end{array}\right.$$
where  $\nu>0$ is the consistency of the fluid, $r'$ is the conjugate exponent of $r$ satisfying $1/r+1/r'=1$, function $f_1$ is the external force and function $g$ is such that its graph defines the lower boundary of the thin film (both functions defined in $\omega$).\\

 {\bf Novelties in the techniques.} We point out that the version of the unfolding method and the restriction operator, introduced in this paper to study the asymptotic behavior of the fluid in the thin porous medium  $\Omega_\ep$, are  powerful tools that could be used to derive lower-dimensional macroscopic laws for different type of (two dimensional or three dimensional) non-Newtonian fluids in a thin porous medium. \\

{\bf Future improvements. }Using the present study as a starting point, various improvements can be proposed. The first one is the generalization of the asymptotic study, which leads to the coupled Darcy--Reynolds equation, to a truly (stationary or non-stationary) nonlinear $3D$ Navier-Stokes system (and not only Stokes system). Another possible way is regarding the boundary conditions. To avoid technical difficulties connected with non-homogeneous boundary conditions for velocity (or pressure in some cases), we have considered a flow with no-slip condition on the exterior boundary of the domain. To derive a more general limit problem, we remark that, with some technical efforts,  this model could be adapted to periodic boundary conditions on the lateral boundaries, to the case of a non-Newtonian fluid with injection as in \cite{MikelicHornung}, or to stress (Neumann) boundary condition on the lateral boundary as in \cite{Fabricius3, Fabricius0, Fabricius}.

\section*{Acknowledgements}
Mar\'ia dedicates this paper to her father, Julio, for all his infinite love. The authors would like to thank the anonymous referees for their nice comments that have allowed us to improve this article.

 \section*{Conflict of interest}
 The authors confirm that there is no conflict of interest to report.

 \section*{Data availability statement}
No new data/code were created or analyzed in this study.

\end{document}